\documentclass[12pt,a4paper]{article}
\usepackage[utf8]{inputenc}
\usepackage[english]{babel}
\usepackage[usenames, dvipsnames]{xcolor}
\usepackage{amssymb}
\usepackage{amsfonts}
\usepackage{amsmath}
\usepackage{amsthm}
\usepackage{epsfig}
\usepackage{dsfont}
\usepackage{enumerate}

\usepackage{diagbox} 

\usepackage{mathtools} 

\usepackage{verbatim}

\usepackage{latexsym}

\usepackage{todonotes}
\usepackage{subfigure}

\allowdisplaybreaks[4]

\theoremstyle{plain}
\newtheorem{theorem}{Theorem}[section]
\newtheorem{lemma}[theorem]{Lemma}
\newtheorem{proposition}[theorem]{Proposition}
\newtheorem{remark}[theorem]{Remark}
\newtheorem{corollary}[theorem]{Corollary}

\theoremstyle{definition}

\newcommand{\rmd}{\, \mathrm d}
\newcommand{\rmdd}{\mathrm d}
\newcommand{\ind}[1]{\mathbf{1}_{#1}}
\newcommand{\ew}{\mathbf{E}}
\newcommand{\var}{\mathbf{Var}}
\newcommand{\cov}{\mathbf{Cov}}
\newcommand{\cor}{\mathbf{Corr}}
\newcommand{\cum}{\mathbf{Cum}}
\newcommand{\F}{\mathcal{F}}
\newcommand{\A}{\mathcal{A}}
\newcommand{\cC}{\mathcal{C}}

\newcommand{\R}{\mathbb{R}}

\newcommand{\N}{\mathbb{N}}
\newcommand{\Z}{\mathbb{Z}}

\newcommand{\nn}{\nonumber}
\newcommand{\tn}[1]{\textup{#1}}

\newcommand{\mb}[1]{\mathbf{#1}}
\newcommand{\X}{\mb{X}}

\newcommand{\e}{\tn{e}}

\newcommand{\bi}{\begin{itemize}}
\newcommand{\ei}{\end{itemize}}

\newcommand{\lbd}{\lambda}

\newcommand{\ds}{\displaystyle}

\newcommand{\wtl}[1]{\widetilde{#1}}

\newcommand{\skp}[2]{\langle #1,#2 \rangle}
\newcommand{\norm}[1]{\left\Vert #1 \right\Vert}

\newcommand{\dcon}{\overset{d}{\longrightarrow}}

\renewcommand{\[}{\left[}
\renewcommand{\]}{\right]}

\begin{document}
\title{On the sample autocovariance of a L\'evy driven moving average process when sampled at a renewal sequence}
\author{Dirk-Philip Brandes\thanks{Ulm University,
Institute of Mathematical Finance, D-89081 Ulm, Germany,
eMail: dirk.brandes@uni-ulm.de, tel.: +49/731/50-23654.} \quad \quad Imma Valentina Curato\thanks{Ulm University,
Institute of Mathematical Finance, D-89081 Ulm, Germany,
eMail: imma.curato@uni-ulm.de, tel.: +49/731/50-23518.}}
\date{\today}
\maketitle
\begin{abstract}
We consider a L\'evy driven continuous time moving average process $X$ sampled at random times which follow a renewal structure independent of $X$. Asymptotic normality of the sample mean, the sample autocovariance, and the sample autocorrelation is established under certain conditions on the kernel and the random times. We compare our results to a classical non-random equidistant sampling method and give an application to parameter estimation of the L\'evy driven Ornstein-Uhlenbeck process. 
\end{abstract}
\noindent{\bf Mathematics subject classification:} 60F05, 60G10, 62D05.  \\

\noindent{\bf Key words:}  central limit theorem, continuous time moving average process, L\'evy process, Ornstein-Uhlenbeck process, renewal sampling, sample autocovariance, sample autocorrelation, sample mean, stationarity.

\section{Introduction}



The present paper analyzes the distributional limit of sample mean and sample autocovariance function of a L\'evy driven moving average process when sampled at a renewal sequence. More precisely, let $X=(X_t)_{t \in \R}$ be a continuous time moving average process of the form
\begin{align}\label{eq1.1}
X_t = \mu + \int_{\R} f(t-s) \rmd L_s \, , \quad t \in \R \, ,
\end{align}
where $L=(L_t)_{t \in \R}$ is a two-sided $\R$-valued L\'evy process, $\mu \in \R$, and $f \colon \R \to \R$ a deterministic function, called \emph{kernel}, for which the integral exists.

Moving average processes as in (\ref{eq1.1}) are the natural continuous time analogue of discrete time moving average processes, whose distributional limits for sample mean and sample autocorrelation have been studied by Brockwell and Davis \cite{brock:06}, Davis and Mikosch \cite{DM:98}, and Hannan \cite{han:76}, and many others.

A popular example of a L\'evy driven moving average process is given by the Ornstein-Uhlenbeck (OU) process used to model the volatility of a financial asset, see \cite{bns:01}, or the intermittency in a turbulence flow, see \cite{BNSCH:07}. The OU process is in fact a tractable mathematical model that can adequately describe the fluctuations of the price and the velocity field on different time scales.



The process $X$ is infinitely divisible, as shown in Rajput and Rosinski \cite{RR:89}, and strictly stationary meaning that its finite dimensional joint distributions are shift-invariant, i.e. for all $n \in \N$ and all $t_1,\dots,t_n \in \R$ it holds
$$\mathcal{L}(X_{t_1+h},\dots,X_{t_n+h}) = \mathcal{L}(X_{t_1},\dots,X_{t_n}) \quad \forall \, h \in \R \, .$$

We study a renewal sampling of the process $X$ in (\ref{eq1.1}). We select a sequence of increasing random times $(T_n)_{n \in \Z}$ such that $T_n \to \infty$ almost surely (abbreviated a.s.). More in detail, we assume that $W=(W_n)_{n \in \Z\setminus\{0\}}$ is an i.i.d. sequence of positive supported random variables independent of the driving L\'evy process $L$ and such that $P(W_1 > 0) > 0$. We then define $(T_n)_{n \in \Z}$ by
\begin{align}\label{eq1.2}
T_0:=0 \quad \text{and} \quad T_n := \begin{cases}
\sum_{i=1}^n W_i \, , &\quad n \in \N \, , \\
-\sum_{i=n}^{-1} W_i \, , &\quad -n \in \N \, ,
\end{cases}
\end{align}
and the sampled process $Y=(Y_n)_{n \in \Z}$ via
\begin{align}\label{eq1.3}
Y_n := X_{T_n} \, , \quad n \in \Z \, .
\end{align}

We are interested in studying the sample moments of the process $Y$. We do this for different reasons. First of all, continuous processes are often used in time series analysis because they can be sampled at non-equidistant points in time and therefore provide a model for non-equidistant data which are often available for statistical inference. 
Secondly, several authors (Cohen and Lindner \cite{CL:13}, Drapatz \cite{dr:17}, and Spangenberg \cite{sp:15}) analyze the asymptotic distribution of the sample mean and sample autocovariance function when $X$ is sampled on a lattice $\{\Delta t \colon t=0,1,2,\dots\}$ for $\Delta > 0$ but results for non-equidistant sampling schemes have not yet been shown.


The central limit theorems presented in the paper generalize the results of Cohen and Lindner \cite{CL:13} at the cost of slightly more restrictive moment conditions. In fact, in the latter, the asymptotic normality of the sample mean and the sample autocovariance function is shown under the assumption of finite second and fourth moments of the driving L\'evy process, respectively, whereas we achieve the asymptotic normality by requiring $\ew(|L_1|^2\log^+|L_1|)$ and $\ew(|L_1|^4 (\log^+|L_1|)^2)$, respectively, to be finite. 

As an application of the developed theory, we present a parameter estimation of the mean reverting parameter of a L\'evy driven OU process
\begin{align}\label{eq*14}
X_t = \int_{-\infty}^t \e^{-a(t-s)} \rmd L_s \, , \quad t \in \R \, ,
\end{align}
sampled at a Poisson process, i.e. a sequence $(T_n)_{n \in \Z}$ where $W$ is a sequence of i.i.d. exponentially distributed random variables. We then compare the efficiency of our estimator with an estimator based on the results of Cohen and Lindner \cite{CL:13} for an equidistant sampling.

The paper is organized as follows. In Section 2 we give some preliminary results regarding strict stationarity of a process sampled at a renewal sequence and the mixing property that it fulfills. Section 3 is concerned with establishing a central limit theorem for the sample mean of a renewal sampled continuous time moving average process as Section 4 does for the sample autocovariance and sample autocorrelation functions. Finally in Section 5, we show the parameter estimation of a L\'evy driven OU process, and Section 6 concludes.




\section{Preliminaries}
\setcounter{equation}{0}

We start with some preliminary results on the continuous time moving average process and its renewal sampled processes needed in the upcoming sections. As a first results, we prove that a strictly stationary process sampled at a renewal sequence inherits the strict stationarity. In particular, this shows that the sampled process (\ref{eq1.3}) is strictly stationary.

Denote with $\overset{d}{=}$ equality in distribution and with $A'$ the transpose of a matrix $A \in \R^{d\times m}$. 

\begin{proposition}\label{prop1.1}
Let $X=(X_t)_{t \in \R}$ be an $\R^d$-valued strictly stationary process $X_t=(X_t^{(1)},\dots,X_t^{(d)})'$ and let $(T_n)_{n \in \Z}$ be a sequence of random times defined as in (\ref{eq1.2}) independent of $X$.
Then the $\R^d$-valued process $Y=(Y_n)_{n \in \Z}$ defined by $Y_n^{(i)}:=X_{T_n}^{(i)}$, $i=1,\dots,d$, is strictly stationary. More generally, the process $((Y_n',T_n-T_{n-1})')_{n \in \Z}$ is strictly stationary.
\end{proposition}

\proof Observe that $((Y_n',T_n-T_{n-1})')_{n \in \Z}$ is strictly stationary if and only if $((Y_n',T_{n+1}-T_{n})')_{n \in \Z}$ is strictly stationary, and the latter implies strict stationarity of $Y$. Hence, it suffices to show that $((Y_n',T_{n+1}-T_{n})')_{n \in \Z}$ is strictly stationary. Let $m \leq n$, $B \in \mathcal{B}(\R^{(d+1)(n-m+1)})$, the Borel-$\sigma$-algebra on $\R^{(d+1)(n-m+1)}$, and denote with $P_Z$ the finite dimensional distribution of a random vector $Z$. Define
$$R_k := T_{k+m}-T_m \, , \quad k =1,\dots, n-m+1 \, .$$
Conditioning and using the strict stationarity of $X$, we obtain
\begin{align*}
& P((Y_m',\dots,Y_n', T_{m+1}-T_m,\dots,T_{n+1}-T_n)' \in B) \\
&\quad= P((X_{T_m}', X_{T_m+R_1}', \dots,X_{T_m+R_{n-m}}',R_1,R_2-R_1,\dots,R_{n-m+1}-R_{n-m})' \in B) \\
&\quad= \int_{\R^{n-m+2}} P((X_u',X_{u+v_1}',\dots,X_{u+v_{n-m}}',v_1,v_2-v_1,\dots,v_{n-m+1}-v_{n-m})' \in B) \\
&\hspace*{8cm} P_{(T_m,R_1,\dots,R_{n-m+1})} (\rmd (u,v_1,\dots,v_{n-m+1})) \\
&\quad= \int_{\R^{n-m+1}} P((X_0',X_{v_1}',\dots,X_{v_{n-m}}',v_1,v_2-v_1,\dots,v_{n-m+1}-v_{n-m})' \in B) \\
&\hspace*{8cm}  P_{(R_1,\dots,R_{n-m+1})} (\rmd (v_1,\dots,v_{n-m+1})) \, ,
\end{align*}
where in the last line we used that the integrand does not depend on $u$. Using that $(R_1,\dots,R_{n-m+1})' \overset{d}{=} (T_1,\dots,T_{n-m+1})'$, the latter is equal to 
\begin{align*}
&=P((Y_0',\dots,Y_{n-m}',T_1,T_2-T_1,\dots,T_{n-m+1}-T_{n-m})' \in B) \, ,
\end{align*}
showing the strict stationarity of $((Y_n',T_{n+1}-T_n)')_{n \in \Z}$. 
\qed

\vspace{0.2cm}

In order to prove central limit theorems, we recall the concept of mixing:
on a probability space $(\Omega,\F,P)$ for any two $\sigma$-algebras $\A,\cC \subset \F$ the following measures of dependence can be defined
\begin{align*}
\alpha(\A,\cC,P) &:= \sup |P(A\cap C)-P(A)P(C)| \, , \quad A \in \A \, , \ C \in \cC \, , \\
\rho(\A,\cC,P) &:= \sup |\cor(f,g)| \, , \quad f \in L^2(\Omega,\A,P) \, , \ g \in L^2(\Omega,\cC,P) \, .
\end{align*}
We say that a strictly stationary sequence of random vectors $Z=(Z_n)_{n \in \Z}$ is
\begin{align*}
\text{\emph{strongly mixing} if} \ \alpha_n &:= \alpha(\A,\cC_n ; P) \to 0 \ \text{ as } n \to \infty \, , \\
\rho\text{-\emph{mixing} if} \ \rho_n &:= \rho(\A,\cC_n ; P) \to 0 \ \text{ as } n \to \infty \, , 
\end{align*}
for the $\sigma$-algebras $\A=\sigma(Z_{0},Z_{-1},\dots)$ and $\cC_n=\sigma(Z_{n},Z_{n+1},\dots)$, cf. Bradley \cite{brad:07-1}.

Recall that a process $X=(X_t)_{t \in \R}$ is called an $m$-dependent process when $(X_t)_{t \leq s}$ and $(X_t)_{t > s+m}$ are independent for each $s$.

\begin{proposition}\label{prop2.1}
Let $X=(X_t)_{t \in \R}$ be for some $m \in \N$ an $\R^d$-valued $m$-dependent strictly stationary process and $Y=(Y_n)_{n \in \Z}$ defined by $Y_n^{(i)}:=X_{T_n}^{(i)}$ with $(T_n)_{n \in \Z}$ as in (\ref{eq1.2}) independent of $X$. Then $Y$ is strongly mixing with exponentially decreasing mixing coefficients $\alpha_n$. More generally, $((Y_n',T_{n}-T_{n-1})')_{n \in \Z}$ is strongly mixing with exponentially decreasing mixing coefficients. 
\end{proposition}

\proof Let $Z_n := (Y_n',T_n-T_{n-1})'$, $n \in \Z$. Let $D_n := \{T_n > m\}$ and $n \geq 1$ so large such that $P(D_n)>0$. First, we show that the two $\sigma$-algebras $\A$ and $\cC_n$ defined as above are independent under the conditional probability measure $P(\cdot | D_n)$. Therefore, let $A \in \A$ and $C_j \in \cC_n$ be of the form $A=\{X_{T_i}' \in B \, , T_i-T_{i-1} \in F \}$ for some $i \leq 0$, $B \in \mathcal{B}(\R^d)$, and $F \in \mathcal{B}(\R)$, and $C_j =\{X_{T_j}' \in B' \, , T_j-T_{j-1} \in F' \}$ for some $j \geq n$, $B' \in \mathcal{B}(\R^d)$, and $F' \in \mathcal{B}(\R)$, respectively. Then, by the Doob-Dynkin lemma and the $m$-dependence of $X$,
\begin{align*}
&P(A\cap C_j | D_n) 
= \frac{1}{P(D_n)} \int_{(m,\infty)} \ew[\ind{A\cap C_j} |T_n=t] P_{T_n}(\rmdd  t) \\
&= \frac{1}{P(D_n)} \int_{(m,\infty)} P(X_{T_i}' \in B , T_i-T_{i-1} \in F , X_{T_j}' \in B', T_j-T_{j-1} \in F' |T_n=t) P_{T_n}(\rmdd  t) \\
&=  P(A | D_n) \frac{1}{P(D_n)}  \int_{(m,\infty)} P( X_{T_j}' \in B',T_j-T_{j-1} \in F'|T_n=t) P_{T_n}(\rmdd  t) \, .
\end{align*}
Observe that $P(A)=P(A|D_n)$ since $X_{T_i}'$ for $i \leq 0$ and $T_n$ are independent.
A calculation like the one above for $B=\R^d$, i.e. $A=\Omega$, gives 
$$P(C_j | D_n) =\frac{1}{P(D_n)} \int_{(m,\infty)} P( X_{T_j}' \in B',T_j-T_{j-1} \in F' | T_n=t) P_{T_n}(\rmdd  t) \\$$
such that all together we obtain
\begin{align}\label{eq2.3}
P(A\cap C_j | D_n) = P(A |D_n) P(C_j |D_n) \quad \text{for } j \geq n \, .
\end{align}
Similarly we can obtain (\ref{eq2.3}) for $A'=\{X_{T_{i_1}}'\in B_1,\dots,X_{T_{i_k}}' \in B_k,T_{i_1}-T_{i_1-1}\in F_1,\dots,T_{i_k}-T_{i_k-1}\in F_k\}$ for $i_1,\dots,i_k \leq 0$, $B_1,\dots,B_k \in \mathcal{B}(\R^d)$, and $F_1,\dots,F_k \in \mathcal{B}(\R)$, and $C_n'=\{X_{T_{j_1}}' \in B_1',\dots,X_{T_{j_l}}' \in B_l',T_{j_1}-T_{j_1-1}\in F_1',\dots,T_{j_l}-T_{j_l-1}\in F_l'\}$ for $j_1,\dots,j_l \geq n$, $B_1',\dots,B_l'\in \mathcal{B}(\R^d)$, and $F_{1}',\dots F_l' \in \mathcal{B}(\R)$.
Observe that sets of the form $A'$ generate the $\sigma$-algebra $\A$ and sets of the form $C_n'$ generate the $\sigma$-algebra $\cC_n$ and both are $\cap$-stable. Thus, we conclude that (\ref{eq2.3}) is true for all $A \in \A$ and $C_n \in \cC_n$. Using measure theoretic induction, and
\begin{align*}
\cov_{P(\cdot|D_n)}(\ind{A},\ind{C_n}) = P(A\cap C_n|D_n) - P(A|D_n)P(C_n|D_n) = 0 \, ,
\end{align*}
we obtain that $\rho(\A,\cC_n,P(\cdot|D_n))=\sup|\cor_{P(\cdot|D_n)}(f,g)| = 0$ where the supremum is taken over all $f \in L^2(\Omega,\A,P(\cdot|D_n))$ and $g \in L^2(\Omega,\cC_n,P(\cdot|D_n))$.

Since $P(D_n) = 1-P(D_n^c)$ and $0=\rho(\A,\cC_n;P(\cdot|D_n)) \leq P(D_n^c)$, it follows, from a remark between Theorem 1.2 and Theorem 1.3 in Bradley \cite{brad:90}, that
$$\alpha_n(\A,\cC_n;P) \leq 4P(D_n^c) = 4P(T_n \leq m) \, .$$

Since $(W_n)_{n \in \Z\setminus\{0\}}$ is supported on $[0,\infty)$ and $P(W_1 > 0 ) > 0$, there exists an $r > 0$ such that $P(W_1 > r )> 0$ and hence $P(W_1+\dots+W_{\lceil m/r \rceil} > m) > 0$, where $\lceil x \rceil$ for $x \in \R$ denotes the smallest integer $k \in \N$ so that $k \geq x$. Denote $q:=1-P(W_1 + \dots + W_{\lceil m/r \rceil} > m) < 1$. Then, as long as $n \leq \lceil m/r \rceil$, we obtain $P(T_n \leq m) \leq q$. For $n > \lceil m/r \rceil$ set $k_n=\lfloor \frac{n}{\lceil m/r \rceil} \rfloor$ for $n \in \N$. Then, by the i.i.d. property of $(W_n)_{n \in \Z\setminus\{0\}}$,
\begin{align}
P(T_n \leq m) \leq  
P(W_1+\dots+W_{\lceil m/r \rceil} \leq m)^{k_n} P(W_{k_n\lceil m/r \rceil +1}+\dots+W_n \leq m ) 
\leq  q^{k_n} \, , \nn
\end{align}
and
\begin{align}\label{eq2.4}
\alpha_n(\A,\cC_n,P) \leq 4q^{k_n} \to 0 \quad \text{as } \ n \to \infty 
\end{align}
showing that $Z$ and hence $Y$ are strongly mixing with exponentially decreasing mixing coefficients.
\qed 

\vspace{0.2cm}

We give some results on when a L\'evy driven continuous time moving average process is well-defined. Moreover, we show under which conditions finiteness of the moments of the process is achieved. 


An $\R$-valued L\'evy processes $L=(L_t)_{t \in \R}$, i.e. a stochastic process with independent and stationary increments, c\`adl\`ag sample paths and $L_0=0$ almost surely, which is continuous in probability, is characterized by its characteristic triplet $(\sigma^2_L, \nu_L, \gamma_L)$ due to the L\'evy-Khintchine formula, i.e.
if $\mu$ denotes the infinitely divisible distribution of $L_1$, then its characteristic function is given by
\begin{align*}
\widehat{\mu}(z) = \exp\[ i \gamma_L z -\frac{1}{2} \sigma^2_L z^2 +\int_{\R} (\e^{izx} - 1 - izx \, \ind{\{|x|\leq 1\}} ) \, \nu_L (\rmdd  x) \] \, , \quad z \in \R \, .
\end{align*}
Here, $\sigma^2_L$ is the Gaussian covariance, $\nu_L$ a measure on $\R$ which satisfies $\nu_L(\{0\}) = 0$ and $\int_{\R} (|x|^2 \land 1) \, \nu_L(\rmdd  x) < \infty$, called the \emph{L\'evy measure}, and $\gamma_L \in \R$.

For a detailed account on L\'evy processes we refer to the book of Sato \cite{sato:13}. 

In the following theorem, we recall the multivariate extension of Theorem 2.7 in Rajput and Rosinski \cite{RR:89}, which characterizes the continuous time moving average process.

\begin{theorem}\label{thm3.1}
Let $L=(L_t)_{t \in \R}$ be a L\'evy process on $\R$ with characteristic triplet $(\gamma_L,\sigma^2_L,\nu_L)$ and $g \colon \R \to \R^d$ be a measurable function. Denote with $D_d:=\{x \colon |x| \leq 1\}$ the unit ball in $\R^d$. Then 

\noindent (a) $g$ is $L$-integrable $($i.e. integrable with respect to the L\'evy process $L)$ as a limit in probability in the sense of Rajput und Rosinski \cite{RR:89} if and only if \\
\vspace{-2mm}

(i) $ \int_\R \norm{g(s)\gamma_L + \int_\R g(s)x (\ind{D_d}(g(s)x) - \ind{D_1}(x)) \, \nu_L(\rmdd x)} \rmd s < \infty$, \\

(ii) $\int_\R \norm{g(s) \sigma^2_L g(s)'} \rmd s < \infty$, and \\

(iii) $\int_\R \int_\R (\norm{g(s)x}^2 \land 1) \, \nu_L(\rmdd x) \rmd s < \infty$. \\

\noindent (b) If $g$ is $L$-integrable, the distribution of $\int_\R g(s) \rmd L_s$ is infinitely divisible with characteristic triplet $(\gamma_{\mathrm{int}},\Sigma_{\mathrm{int}},\nu_{\mathrm{int}})$ given by
\begin{align*}
\gamma_{\mathrm{int}} &= \int_\R g(s) \gamma_L + \int_{\R} g(s)x (\ind{D_d}(g(s)x)-\ind{D_1}(x)) \, \nu_L(\rmdd x) \rmd s \, ,\\
\Sigma_{\mathrm{int}} &= \sigma^2_L \int_\R g(s) g(s)' \rmd s \, , \quad \text{and}\\
\nu_{\mathrm{int}} (B) &= \int_\R \int_\R \ind{B}(g(s)x) \, \nu_L(\rmdd x) \rmd s \quad \text{for all Borel sets } B \subset \R^d\setminus\{0\}.
\end{align*}
\end{theorem}

\begin{corollary}\label{cor3.2}
By a simple calculation, we deduce that, if $L$ has expectation zero and finite second moment and $g \in L^2(\R^d)$, then the conditions (i), (ii), and (iii) of Theorem \ref{thm3.1} (a) are satisfied, $g$ is $L$-integrable, and $\int_\R g(s) \rmd L_s$ is infinitely divisible with characteristic triplet $(\gamma_{\mathrm{int}},\Sigma_{\mathrm{int}},\nu_{\mathrm{int}})$ as given in Theorem \ref{thm3.1} (b).
\end{corollary}

The next lemma shows, for a process $X$ as in (\ref{eq1.1}), finiteness of the $\log$-moment under certain conditions on the L\'evy process $L$ and the kernel function $f$.

We use the notation $\log^+(x):=\log(\max\{1,x\})$.

\begin{lemma}\label{lem2.2}
Let $X=(X_t)_{t \in \R}$ be defined by $X_t := \mu + \int_\R f(t-s) \rmd L_s \, ,$ where $f \in L^2(\R)$ and $L=(L_t)_{t \in \R}$ is a one-dimensional L\'evy process with zero mean. \\
%
\vspace*{-0.4cm}

\noindent(a) If $\ew(|L_1|^2\log^+|L_1|)<\infty$, and $\int_\R |f(s)|^2 \log^+|f(s)| \rmd s < \infty$,
then $\ew(|X_t|^2\log^+|X_t|)<\infty$ for all $t \in \R$.\\
\vspace*{-0.4cm}

\noindent(b) If $\ew(|L_1|^4(\log^+|L_1|)^2)<\infty$, $f \in L^4(\R)$, and $\int_\R |f(s)|^4 (\log^+|f(s)|)^2 \rmd s < \infty $,
then $\ew(X_t^4(\log^+|X_t|)^2) < \infty$ for all $t \in \R$ and, for $h \in \N$, $\ew(|X_tX_{t+h}|^2\log^+|X_tX_{t+h}|)<\infty$ for all $t \in \R$.
\end{lemma}

\proof W.l.o.g. $\mu=0$. It is enough to show the assertions for $Z=\int_\R f(s) \rmd L_s$, for which $Z\overset{d}{=}\int_\R f(-s) \rmd L_s =X_0$. By the strict stationarity of $X$, we obtain the result.

(a) By Corollary \ref{cor3.2}, $Z$ is infinitely divisible with triplet $(\gamma_Z,\sigma^2_Z,\nu_Z)$ given by Theorem \ref{thm3.1} (b). By Corollary 25.8 of Sato \cite{sato:13}, we know that $\ew(Z^2\log^+|Z|)<\infty$, if $\int_{|x|>1} |x|^2 \log^+|x| \, \nu_Z(\rmdd  x) < \infty$. To see that this is indeed true, 
observe that $\log^+|ab| \leq \log^+|a| + \log^+|b|$ for $a,b \in \R$. Hence, by the given assymptions,
\begin{align*}
&\int_{|x|>1} |x|^2 \log^+|x| \, \nu_Z(\rmdd  x) =\int_\R \int_\R |f(s)x|^2 \log^+|f(s)x| \, \nu_L(\rmdd  x) \rmd s \\
&\quad\leq \int_\R |f(s)|^2\log^+|f(s)| \rmd s \int_\R |x|^2 \, \nu_L(\rmdd  x) + \int_\R |f(s)|^2 \rmd s \int_\R |x|^2 \log^+|x| \, \nu_L(\rmdd  x) < \infty \, .
\end{align*}

(c) Observe that $\ew(|L_1|^4(\log^+|L_1|)^2)<\infty$ is equivalent to $\int_{|x|>1} |x|^4 (\log^+|x|)^2 \, \nu_L (\rmdd x) < \infty$ and that $\ew(|L_1|^4 (\log^+|L_1|)^2) < \infty$ implies $\ew|L_1|^4 < \infty$.
Henceforth we obtain
\begin{align*}
&\int_{|x|>1} |x|^4 (\log^+|x|)^2 \, \nu_Z(\rmdd  x) =\int_\R \int_\R |f(s)x|^4 (\log^+|f(s)x|)^2 \, \nu_L(\rmdd  x) \rmd s \\
& \ \leq 2 \int_\R \int_\R |f(s)x|^4 ((\log^+|f(s)|)^2 + (\log^+|x|)^2 ) \, \nu_L(\rmdd  x) \rmd s \\
&\ \leq \int_\R |f(s)|^4(\log^+|f(s)|)^2 \rmd s \int_\R |x|^4 \, \nu_L(\rmdd  x) + \int_\R |f(s)|^4 \rmd s \int_\R |x|^4 (\log^+|x|)^2 \, \nu_L(\rmdd  x) < \infty \, .
\end{align*}

This gives $\ew(|X_t|^4(\log^+|X_t|)^2)<\infty$ for all $t \in \R$. 
Using the strict stationarity of $X$ and the Cauchy-Schwarz inequality yields
\begin{align*}
\ew(|X_tX_{t+h}|^2\log^+|X_tX_{t+h}|) 
\leq 2 \ew(|X_0|^4 (\log^+|X_0|)^2) \ew|X_0|^4 < \infty 
\end{align*}
which gives the result.
\qed 

\begin{remark}
Throughout the paper, we assume that the L\'evy process $L=(L_t)_{t \in \R}$ has zero mean. This assumption can be dropped in many cases. For example, if $f \in L^1(\R)\cap L^2(\R)$, we define with $L_t'=L_t - t\ew(L_1)$, $t \in \R$, another L\'evy process with mean zero and the same variance such that
$$X_t = \mu + \ew(L_1)\int_{\R}f(s)\rmd s + \int_\R f(t-s) \rmd L_s' \, , \quad t \in \R \, ,$$
and $X_t$ has mean $\mu + \ew(L_1)\int_{\R}f(s)\rmd s$.
\end{remark}

\section{Sample Mean}
\setcounter{equation}{0}

In this section, we show the asymptotic normality of the sample mean
\begin{align}\label{eq*31}
\overline{Y}_n := \sum_{k=1}^n Y_k = \sum_{k=1}^n X_{T_k} \, , \quad n \in \N \, ,
\end{align}
where $X=(X_t)_{t \in \R}$ and $Y=(Y_n)_{n \in \Z}$ are given in (\ref{eq1.1}) and (\ref{eq1.3}), respectively.

To do so, we consider a certain truncated continuous time moving average process. Therefore, let $f_m\colon\R \to \R, s \mapsto f(s)\ind{[-m/2,m/2]}$ be a kernel function with compact support, and $X^{(m)}=(X^{(m)}_t)_{t \in \R}$ be defined by
\begin{align}\label{eq2.1}
X_t^{(m)} := \mu + \int_{\R} f_m(t-s) \rmd L_s =  \mu + \int_{\R} f(t-s) \ind{[-m/2,m/2]}(t-s) \rmd L_s\, , \quad t \in \R \, ,
\end{align}
where $L=(L_t)_{t \in \R}$ is a L\'evy process with zero mean and $\ew|L_1|^2<\infty$, $\mu \in \R$, and $f \in L^2(\R)$. Then the process $X^{(m)}=(X_t^{(m)})_{t \in \R}$ is an $m$-dependent process. Moreover, $X^{(m)}$ is strictly stationary and, by Proposition \ref{prop1.1}, so is the sequence $Y^{(m)}=(Y^{(m)}_n)_{n \in \Z}$ defined by 
\begin{align}\label{eq2.2}
Y^{(m)}_n:=X^{(m)}_{T_n} \, , 
\end{align}
where $(T_n)_{n \in \Z}$ is defined as in (\ref{eq1.2}) independent of $X$.


The following proposition states a result on the convergence of the covariances of $Y^{(m)}$ towards the ones of $Y$. By Proposition \ref{prop1.1}, the process $Y$ is strictly stationary.

\begin{proposition}\label{prop2.5}
Let $X$ be defined by (\ref{eq1.1}), $X^{(m)}$ by (\ref{eq2.1}), the processes $Y$ and $Y^{(m)}$ by (\ref{eq1.3}) and (\ref{eq2.2}), respectively, with $(T_n)_{n \in \Z}$ as in (\ref{eq1.2}) and assume that $\mu=0$. Then
\begin{align}\label{eq2.5}
\ew(|Y_kY_l - Y^{(m)}_kY^{(m)}_l|) \to 0 \quad \text{as } \ m \to \infty  \quad \text{for } \ k,l \in \Z \, .
\end{align}
Further, it holds
\begin{align}\label{eq2.6}
\ew(Y_kY_l) = \ew(L_1^2) \int_\R \ew(f(u) f(T_{|l-k|}+ u)) \rmd u \quad \text{for } \  k,l \in \Z \, , 
\end{align}
and similar for $\ew(Y^{(m)}_kY^{(m)}_l)$ with $f$ replaced by $f_m=f\ind{[-m/2,m/2]}$.
\end{proposition}

\proof 
We denote with $\sigma(T)$ the $\sigma$-algebra generated by some random variable $T$. Clearly, $|f_m(u)| \leq |f(u)|$ for all $u \in \R$. Then, by conditioning on $T_k$ and Fubini's theorem,
\begin{align}\label{eq2.8}
\ew|Y_k^{(m)}|^2 
&=\int_\R \ew\bigg( \bigg(\int_\R f_m(t-u)\rmd L_u\bigg)^2 \bigg| T_k=t\bigg) P_{T_k}(\rmdd  t) \nn \\
&=\int_\R \ew(L_1^2) \int_\R (f_m(t-u))^2\rmd u \,  P_{T_k}(\rmdd  t) \leq \ew(L_1^2) \int_\R f(u)^2 \rmd u < \infty \, .
\end{align}

Further, observe that
\begin{align}\label{eq2.9}
\ew|Y_k-Y_k^{(m)}|^2 
=\ew \left[ \ew \left( \left. \left( \int_{\R\setminus [T_k-m/2,T_k+m/2]}f(T_k-u) \rmd L_u \right)^2 \right| \sigma(T_k) \right) \right] =: \ew(\mathrm{I}) \, .
\end{align}

Define
$$\varphi^{(m)} (t) := \ew \left( \left. \left( \int_{\R\setminus [t-m/2,t+m/2]}f(T_k-u) \rmd L_u \right)^2 \right| T_k=t \right) \, ,$$
then obviously $\varphi^{(m)}\circ T_k=\mathrm{I}$. But, since $L$ is independent of $(T_n)_{n \in \Z}$, it holds
\begin{align*}
\varphi^{(m)} (t) 
= \ew(L_1^2) \int_{\R\setminus[t-m/2,t+m/2]} f(t-u)^2 \rmd u \to 0 \quad \text{as } \ m \to \infty \, ,
\end{align*}
since $f \in L^2(\R)$.
Hence $\varphi^{(m)}(T_k(\omega)) \to 0$ as $m \to \infty$ for all $k \in \Z$ and all $\omega \in \Omega$.

Define $\varphi(t):=
\ew(L_1^2)\int_\R f(t-u)^2 \rmd u$.
Then $\ew(\varphi \circ T_k) = \ew(Y_k^2) < \infty$ such that, since $|\varphi^{(m)} \circ T_k| \leq |\varphi\circ T_k|$, we obtain by the dominated convergence theorem for (\ref{eq2.9})
\begin{align}\label{eq2.10}
\ew|Y_k-Y_k^{(m)}|^2 = \ew(\mathrm{I}) = \ew(\varphi^{(m)} \circ T_k)   \to 0 \quad \text{as } m \to \infty \, .
\end{align}

Henceforth, by (\ref{eq2.8}), (\ref{eq2.10}), and the Cauchy-Schwarz inequality,
\begin{align*}
\ew |Y_kY_l - Y^{(m)}_kY^{(m)}_l| &= \ew|Y_kY_l - Y^{(m)}_kY^{(m)}_l + Y^{(m)}_kY_l -  Y^{(m)}_kY_l| \\
& \leq \sqrt{\ew| Y_k^{(m)}|^2} \sqrt{ \ew|Y_l-Y^{(m)}_l|^2} + \sqrt{\ew|Y_l|^2}\sqrt{\ew|Y_k-Y_k^{(m)}|^2}  \to 0 
\end{align*}
for $m \to \infty$, i.e. (\ref{eq2.5}).

For the last statement (\ref{eq2.6}), let w.l.o.g. $k,l \in \N_0$ and $k\leq l$. Then 
\begin{align*}
\ew &(Y_kY_l) 
 = \int_\R \ew\left[ \left. \int_\R f(t-u) \rmd L_u \int_\R f\bigg(t+\sum_{i=k+1}^lW_i-u\bigg) \rmd L_u \right| T_k = t \right] \, P_{T_k}(\rmdd  t) \\
& = \int_{[0,\infty)}\int_{[0,\infty)}  \int_\R f(t-u)  f(t+s-u) \rmd u \, \ew(L_1^2) \,  P_{\sum_{i=k+1}^lW_i} (\rmdd s)P_{T_k}(\rmdd  t) \\
&= \ew(L_1^2)\int_\R \ew(f(v)  f(T_{l-k}+v)) \rmd u \, .
\end{align*}
such that the statement follows.
\qed 

\vspace{0.2cm}

Now we are in the position to prove the asymptotic normality of $\overline{Y}_n$ in (\ref{eq*31}). We denote with $\dcon$ convergence in distribution.

\begin{theorem}\label{thm2.6}
Let $X$ be defined as in (\ref{eq1.1}) such that $\mu \in \R$, $L$ has expectation zero and $\ew(|L_1|^{2}\log^+|L_1|)<\infty$, $f \in L^2(\R)$, and $\int_\R |f(s)|^2\log^+|f(s)|\rmd s < \infty$. 
Let $Y$ be defined by (\ref{eq1.3}) with $(T_n)_{n \in \Z}$ as in (\ref{eq1.2}) 
independent of $L$.
Assume that
\begin{align}\label{eq2.11}
\int_\R |f(u)| \sum_{k=1}^\infty \ew|f(T_k+u)| \rmd u < \infty \, .
\end{align}
Then

\noindent (a) $\sigma_{\overline{Y}}^2:=\sum_{k\in\Z} \cov(Y_0,Y_k)$ exists in $[0,\infty)$, is absolutely convergent, and 
\begin{align}\label{eq2.12}
\sigma_{\overline{Y}}^2 = \ew(L_1^2) \sum_{k\in\Z} \int_\R f(u) \ew(f(T_k+u)) \rmd u \, .
\end{align}
(b) $\ds \sqrt{n} \, (\overline{Y}_{n}-\mu) \dcon N(0,\sigma_{\overline{Y}}^2)$ as $n \to \infty$.
\end{theorem}

\proof Define $\wtl{X}_t=X_t-\mu$ such that with $\wtl{Y}_n=Y_{n}-\mu$ due to the strict stationarity of $(Y_n)_{n \in \Z}$ and since $Y_0=X_0$, we obtain a sequence with expectation zero. Hence, w.l.o.g. $\mu=0$.

(a) Observe that $\ew|Y_0^2| = \ew|X_0|^2 < \infty$ since $f \in L^2(\R)$ and $L$ has finite second moment. Further, by (\ref{eq2.6}) and (\ref{eq2.11}) together with the dominated convergence theorem 
\begin{align}\label{eq2.13}
\sum_{k\in\Z} |\ew( Y_0Y_k) | 
\leq \ew(L_1^2) \int_{\R} |f(u)| \sum_{k\in\Z} \ew|f(T_k+u)| \rmd u < \infty \, . 
\end{align}
This gives the absolute summability of $\sigma_{\overline{Y}}^2$ and a similar calculation without the modulus gives (\ref{eq2.12}).

(b) Let $\overline{Y}^{(m)}_n=\frac{1}{n} \sum_{k=1}^n Y^{(m)}_{k}$, where the sequence $(Y^{(m)}_n)_{n \in \Z}$ is defined as in (\ref{eq2.2})  via the $m$-dependent process $(X^{(m)}_t)_{t \in \R}$ given in (\ref{eq2.1}). The assumptions on $L$ and $f$ imply, by Lemma \ref{lem2.2} (b), $\ew(|Y^{(m)}_0|^{2}\log^+|Y_0^{(m)}|)=\ew(|X^{(m)}_0|^{2}\log^+|X_0^{(m)}|)<\infty$. 

Observe that $(Y^{(m)}_n)_{n \in \Z}$ fulfils the assumptions of Proposition \ref{prop2.1} and it is therefore strongly mixing with exponentially decreasing mixing coefficient $\alpha_n^{Y^{(m)}}$. Hence, due to (\ref{eq2.4}), $\alpha_n^{Y^{(m)}} = O(\e^{k_n\log{q}})$ as $n \to \infty$.
Hence, by Corollary 10.20 (c) in Bradley \cite{brad:07-1}, we obtain 
\begin{align}\label{eq2.14}
\sqrt{n} \, \overline{Y}^{(m)}_n \dcon Z^{(m)} \quad \text{with} \quad Z^{(m)} \overset{d}{=} N(0,\sigma_{\overline{Y}^{(m)}}^2) \, ,
\end{align}
where $\sigma_{\overline{Y}^{(m)}}^2:=\sum_{k\in\Z} \cov(Y^{(m)}_0,Y^{(m)}_k)$
exists in $[0,\infty)$ and is absolutely convergent.

By Proposition \ref{prop2.5}, we have that $\ew(Y_0^{(m)}Y_{k}^{(m)}) \to \ew(Y_0Y_k)$ as $m \to \infty$ and, since 
$$\sum_{k=1}^\infty |\ew(Y_0^{(m)}Y_k^{(m)})| \leq \ew(L_1^2) \int_\R |f(u)| \sum_{k=1}^\infty \ew |f(T_k+u)| \rmd u < \infty \, ,$$
by (\ref{eq2.13}) and $|f_m| \leq |f|$, we conclude with the dominated convergence theorem that $\lim_{m \to \infty} \sigma_{\overline{Y}^{(m)}}^2 = \sigma_{\overline{Y}}^2$. Hence, 
\begin{align}\label{eq2.15}
Z^{(m)} \dcon Z\quad \text{as } \ m \to \infty \quad \text{with} \quad Z \overset{d}{=} N(0,\sigma_{\overline{Y}}^2) \, .
\end{align}

Define for $k \in \Z$, $Y_k^{f-f_m} := \int_{\R \setminus[T_k-m/2,T_k+m/2]} f(T_k-u) \rmd L_u$.
Then $(Y_n^{f-f_m})_{n \in \Z}$ is strictly stationary, by Proposition \ref{prop1.1}. Further, by Cauchy-Schwarz's inequality,
\begin{align*}
\ew &(Y_0^{f-f_m}Y_k^{f-f_m}) \\
&\leq \bigg( \ew \bigg(\int_{\R\setminus[-m/2,m/2]} f(-u) \rmd L_u \bigg)^2  \ew \bigg(\int_{\R\setminus[T_k-m/2,T_k+m/2]} f(T_k-u) \rmd L_u \bigg)^2 \bigg)^{1/2} \to 0
\end{align*}
as $m \to \infty$ since $f \in L^2(\R)$. Since, by (\ref{eq2.13}),
$$\sum_{k\in\Z} |\ew(Y_0^{f-f_m}Y_k^{f-f_m})| \leq \ew(L_1^2) \int_\R |f(u)| \sum_{k\in\Z} \ew|f(T_k+u)| \rmd u \quad \forall \, m \in \N \, ,$$
the dominated convergence theorem yields 
$\lim_{m\to \infty}\sum_{k\in\Z} |\ew(Y_0^{f-f_m}Y_k^{f-f_m})| = 0$.

Hence, by Theorem 7.1.1 in Brockwell and Davis \cite{brock:06},
\begin{align*}
\lim_{m\to \infty} \lim_{n \to \infty} \var(n^{1/2} (\overline{Y}_n -\overline{Y}^{(m)}_n) ) &= \lim_{m\to \infty} \lim_{n \to \infty} n \var \bigg( \frac{1}{n} \sum_{k=1}^n Y_k^{f-f_m}\bigg) \\
&=\lim_{m \to \infty} \sum_{k\in\Z} \ew(Y_0^{f-f_m} Y_k^{f-f_m}) = 0 \, .
\end{align*}

Chebychef's inequality then gives
$$\lim_{m \to \infty} \limsup_{n \to \infty} P(n^{1/2}|\overline{Y}_n -\overline{Y}^{(m)}_n| > \varepsilon) = 0 \quad \forall \, \varepsilon > 0 \, .$$
Together with (\ref{eq2.14}) and (\ref{eq2.15}), the claim follows by a variant of Slutsky's Lemma, cf. Proposition 6.3.9 of Brockwell and Davis \cite{brock:06}.
\qed

\begin{remark}
When $(T_n)_{n \in \Z}$ is deterministic, i.e. $T_n = \Delta n$ for $n \in \Z$ and some $\Delta > 0$, Cohen and Lindner \cite[Theorem 2.1]{CL:13} established the asymptotic normality of the sample mean under the conditions $\ew(L_1^2)<\infty$, $\ew(L_1)=0$, $f \in L^2(\R)$, and 
\begin{align}\label{CL*1}
\bigg( u \mapsto \sum_{j=-\infty}^\infty |f(u+\Delta j)|\bigg) \in L^2([0,\Delta]) \, .
\end{align}
Observe that (\ref{CL*1}) implies (\ref{eq2.11}) since
\begin{align*}
\int_\R |f(u)| \sum_{k=1}^\infty |f(\Delta k + u)| \rmd u 
\leq \int_0^\Delta \bigg(\sum_{j=-\infty}^\infty |f(u+\Delta j)|\bigg)^2 \rmd u \, . 
\end{align*}
So, Theorem \ref{thm2.6} generalizes Theorem 2.1 in \cite{CL:13} to the case of a renewal sampling sequence $(T_n)_{n \in \Z}$, albeit at the cost of the slightly more restrictive conditions $\ew(|L_1|^2\log^+|L_1|)<\infty$ and $\int_\R |f(s)|^2 \log^+|f(s)| \rmd s < \infty$.
\end{remark}

\begin{remark}\label{rem3.5}
Condition (\ref{eq2.11}) is satisfied, for example, for $|f(u)|\leq K(|u|^{-\alpha} \land 1)$ with $\alpha > 1$ and $K > 0$.

To see this, observe that for some $C_\alpha$ 
\begin{align}\label{eq3.14}
\int_{\R} |f(u)| |f(t+u)| \rmd u \leq C_\alpha (|t|^{-\alpha} \land 1) \, .
\end{align}
Hence, 
\begin{align}\label{eq3.15}
\int_\R |f(u)| \sum_{k=1}^\infty \ew|f(T_k+u)| \rmd u 
&\leq C_\alpha \sum_{k=1}^\infty P(T_k \leq 1) +C_{\alpha}\sum_{k=1}^\infty \ew(T_k^{-\alpha} \ind{\{T_k > 1\}}) \, .
\end{align}

The first sum in (\ref{eq3.15}) converges since, as it was shown in the proof of Proposition \ref{prop2.1}, for all $m \in \N$, and $q:=1-P(W_1 + \dots + W_{\lceil m/r\rceil} > m) <1$, if $k \leq \lceil m/r\rceil$, it holds $P(T_k \leq m) \leq q$. Otherwise for $k > \lceil m/r\rceil$, set $l_k(m)=\lfloor \frac{k}{ \lceil m/r\rceil} \rfloor$ such that $P(T_k \leq m) \leq q^{l_k(m)}$.

For establishing the convergence of the second sum, observe that
\begin{align}\label{eq*3.15}
\sum_{k=1}^\infty \ew(T_k^{-\alpha} &\ind{\{T_k > 1\}}) \leq \sum_{k=1}^\infty \ew(T_k^{-\alpha} \ind{\{T_k^{-\alpha} \leq 1\}}) = \sum_{k=1}^\infty \int_{0}^{1} P( T_k^{-\alpha} \ind{\{T_k^{-\alpha} \leq 1\}}> t) \rmd t \nn \\
&=\sum_{k=1}^\infty\int_{0}^{1} P( 1 \leq T_k < t^{-1/\alpha}) \rmd t  \leq \int_{1}^\infty \alpha v^{-\alpha-1}\sum_{k=1}^\infty P(T_k\leq v) \rmd v \, .
\end{align}

Since $P(T_k \leq v) \leq P(T_k \leq \lceil v \rceil) \leq q^{l_k(\lceil v \rceil)}$ and 
$$l_k(\lceil v \rceil) = \left\lfloor \frac{k}{ \lceil \lceil v \rceil/r\rceil} \right\rfloor \geq \frac{k}{ \lceil \lceil v \rceil/r\rceil} -1 \geq C \frac{k}{v/r} -1 \quad \forall \, k \in \N \, , \ v \geq 1 \, ,$$
for some $C > 0$, we obtain with $\tilde q := q^C$ that
\begin{align}\label{eq*3.16}
\sum_{k=1}^\infty P(T_k\leq v) \leq q^{-1} \sum_{k=1}^\infty q^{C \frac{k}{v/r}} \leq \frac{q^{-1}}{1-\tilde{q}^{r/v}} \, .
\end{align}
Let $\tilde q(x) = \tilde q^{x}$ for $x \geq 0$. Then for $x \in [0,1]$, by the mean value theorem, there exists $\xi \in(0,1)$ such that
$$1-\tilde q^x = \tilde q(0)-\tilde q(x) =  (-x) \tilde q'(\xi) = (-x) \tilde q^{\xi} \log(\tilde q) \geq x \tilde q |\log(\tilde q)| $$
since $\tilde q \in (0,1)$. This yields for $v > r$ that $\frac{1}{1-\tilde q^{r/v}} \leq \frac{v}{r \tilde q |\log (\tilde q)|}$. This together with (\ref{eq*3.16}) gives for (\ref{eq*3.15})
\begin{align*}
\int_{1}^\infty \alpha v^{-\alpha-1}\sum_{k=1}^\infty P(T_k\leq v) \rmd v \leq \int_{1}^\infty \alpha v^{-\alpha} \frac{q^{-1}}{r \tilde q |\log (\tilde q)|} \rmd v < \infty \, ,
\end{align*}
since $\alpha > 1$.
\qed
\end{remark}

\section{Sample Autocovariance}
\setcounter{equation}{0}

In this section, we present a multivariate central limit theorem for the autocovariance and autocorrelation functions when the process is sampled at a renewal sequence. We start by considering the strictly stationary, mean zero process
\begin{align}\label{eq4.1}
X_t = \int_\R f(t-s) \rmd L_s \, , \quad t \in \R \, .
\end{align}

As in the previous section, let $(T_n)_{n \in \Z}$ be a sequence of random times defined by (\ref{eq1.2}), and the sampled process $Y_n = X_{T_n}$ for $n \in \Z$. Recall that for a mean zero process, 
\begin{align}\label{eq4.2}
\gamma_n^*(h) = \frac{1}{n} \sum_{k=1}^n Y_k Y_{k+h} \, , \quad h \in \N_0 \, ,
\end{align} 
is a natural estimator for the autocovariance function. 

In analogy to the proof of the asymptotic normality of the sample mean, we first show that, for the truncated sequence $Y^{(m)}=(Y_n^{(m)})_{n \in \Z}$ as in (\ref{eq2.2}) and (\ref{eq2.1}) with $\mu = 0$, the asymptotic normality can be proved.
%
At this point, we need a series of lemmas that allows us to compute the $4^{th}$-order cumulants of the processes $X$ and $Y$.

Part (a) in the lemma below generalizes expression (3.5) in Cohen and Lindner \cite{CL:13} to non-lattice times and presents a different and quicker proof.

\begin{lemma}\label{lem4.3}
Let $f \in L^2(\R)\cap L^4(\R)$, and $L=(L_t)_{t \in \R}$ be a L\'evy process with expectation zero and finite fourth moment. Denote $\sigma^2:=\ew(L_1^2)$, $\eta:=\sigma^{-4} \ew(L_1^4)$, and $f_m:=f\ind{[-m/2,m/2]}$. Then the following statements hold:\\
\noindent(a) For $X_t := \int_\R f(t-u) \rmd L_u$, we have for all $r,s,t,v \in \R$
\begin{align}\label{eq3.3}
\ew(X_r X_s X_t X_v) &=(\eta-3)\sigma^4 \int_\R f(u+r)f(u+s)f(u+t)f(u+v) \rmd u \nn \\
&+ \ew(X_r X_s)\ew(X_t X_v)  +\ew(X_rX_t) \ew(X_sX_v) + \ew(X_rX_v)\ew(X_sX_t) \, .
\end{align}
\noindent(b) Let additionally $X^{(m)}_t := \int_\R f_m(t-u) \rmd L_u$, then we have for all $r,s,t,v \in \R$
\begin{align*}
\ew(X_r X_s X_t^{(m)} X_v^{(m)}) &=(\eta-3)\sigma^4 \int_\R f(u+r)f(u+s)f_m(u+t)f_m(u+v) \rmd u \\
&\hspace*{-2cm}+ \ew(X_r X_s)\ew(X_t^{(m)} X_v^{(m)})  +\ew(X_rX_t^{(m)}) \ew(X_sX_v^{(m)}) + \ew(X_rX_v^{(m)})\ew(X_sX_t^{(m)}) \, .
\end{align*}
\end{lemma}

\proof 
Since $X_r$, $X_s$, $X_t$, and $X_v$ all have expectation zero, the $4^{th}$ order joint cumulant $\cum(\mathbf{X})$ of $\mathbf{X}:=(X_r,X_s,X_t,X_v)$ is given by 
\begin{align}\label{cum1}
\cum(\mathbf{X}) &= \ew(X_rX_sX_tX_v) - \ew(X_rX_s)\ew(X_tX_v) \nn \\
&-\ew(X_rX_t)\ew(X_sX_v)-\ew(X_rX_v)\ew(X_sX_t) \, ,
\end{align}
see Proposition 4.2.2 in Giraitis et al. \cite{GKS:12}. On the other hand,
$$\cum(\mathbf{X}) = \frac{\partial^4}{\partial u_1 \partial u_2 \partial u_3 \partial u_4} \log \ew(\e^{i \skp{u}{\X}}) \bigg|_{u_1=u_2=u_3=u_4=0} \, ,$$
cf. Definition 4.2.1 of Giraitis et al. \cite{GKS:12}.
Let $g(u)=(f(r-u),f(s-u),f(t-u),f(v-u))$, then it holds
$\X=\int_\R g(u) \rmd L_u$,
and, since $f \in L^2(\R)\cap L^4(\R)$, we obtain $g \in L^2(\R^4) \cap L^4(\R^4)$ which yields by Corollary \ref{cor3.2} that $\X$ is infinitely divisible with characteristic triplet $(\gamma_{\mathbf{X}},\Sigma_{\mathbf{X}},\nu_{\mathbf{X}})$,
by Theorem \ref{thm3.1} (b). This and the well-known fact $\int_\R x^4 \nu_L(\rmdd x) = (\eta-3)\sigma^4$ yield
\begin{align*}
\cum(\X) &= \int_{\R^4} x_1x_2x_3x_4 \, \nu_{\mathrm{int}} (\rmdd x) 
	=(\eta-3)\sigma^4 \int_\R  f(r-u)f(s-u)f(t-u)f(v-u)  \rmd s \, .
\end{align*}
The latter together with (\ref{cum1}) yields (a). 

For (b) just observe that also $f_m \in L^2(\R)\cap L^4(\R)$ such that with $h(u)=(f(r-u),f(s-u),f_m(t-u),f_m(v-u))$ we obtain that also 
$\mathbf{Z}=\int_\R h(u) \rmd L_u$
is infinitely divisible by Corollary \ref{cor3.2}. Similar argumentations as above give the demanded result.
\qed

\vspace{0.2cm}

%

In the following lemma we give a similar expression as (\ref{eq3.3}) when the deterministic times $r,s,t$, and $v$ are replaced by random times. 

\begin{lemma}\label{lem4.5}
Let $L=(L_t)_{t \in \R}$ be a L\'evy process with expectation zero and finite fourth moment, $X$ be defined by (\ref{eq4.1}), with $f \in L^2(\R)\cap L^4(\R)$, and the processes $Y$ is defined by (\ref{eq1.3}) with $(T_n)_{n \in \Z}$ as in (\ref{eq1.2}). Denote $\sigma^2:=\ew(L_1^2)$ and $\eta:=\sigma^{-4}\ew(L_1^4)$, and let $l,m,n \in \Z$. Let $F(s,t):=\int_\R f(u+s)f(u+t) \rmd u$, then \\ 
\vspace*{-0.3cm}

\noindent(a) \vspace*{-10mm}
\begin{align*}
\ew(Y_0Y_lY_mY_n) &= (\eta-3)\sigma^4 \int_\R f(u) \ew(f(u+T_{l})f(u+T_{m})f(u+T_{n})) \rmd u \\
&+ \sigma^4 \ew(F(0,T_{l}))F(T_{m},T_{n})) +  \sigma ^4\ew(F(0,T_{m})F(T_{l},T_{n})) \\
&+ \sigma^4 \ew(F(0,T_{n})F(T_{l},T_{m})) \, .
\end{align*}
\noindent(b) If $0 \leq l \leq m \leq n$, then $\ew(F(0,T_{l}))F(T_{m},T_{n}))=\ew(F(0,T_{l}))\ew(F(0,T_{n-m}))$. 
\end{lemma}

\proof (a) Due to the definition of $(T_n)_{n \in \Z}$ it follows that $T_l \leq T_m \leq T_n$. Conditioning on the random times yields, by the independence of $L$ and $W$ and Lemma \ref{lem4.3} (a),
\begin{align*}
\ew(Y_0Y_lY_mY_n) 
&=\int_{[0,\infty)^3} \ew\big[X_{0}X_{s}X_{t}X_{v} | (T_l, T_m, T_n)'=(s,t,v)'  \big] P_{(T_l,T_m,T_n)} (\rmdd (s,t,v)) \\
&=\int_{[0,\infty)^3} \ew\big(X_{0}X_{s}X_{t}X_{v} ) P_{(T_l,T_m,T_n)} (\rmdd (s,t,v)) 
=: \mathrm{A}+\mathrm{B}+\mathrm{C}+\mathrm{D} \, ,
\end{align*}
where A, B, C, and D correspond to the parts arising from the decomposition in (\ref{eq3.3}). Then, by Fubini's theorem,  
\begin{align*}
\mathrm{A} & = (\eta-3)\sigma^4 \int_{[0,\infty)^3} \int_\R f(u)f(u+s)f(u+t)f(u+v) \rmd u \,  P_{(T_l,T_m,T_n)} (\rmdd (s,t,v)) \\
& =(\eta-3)\sigma^4  \int_\R f(u) \, \ew\big[f(u+T_{l})f(u+T_{m})f(u+T_{n})\big] \rmd u \, .
\end{align*}
Since 
$\ew(X_sX_t) = \sigma^2 \int_\R f(u+s)f(u+t) \rmd u$ for all $s,t \in \R$, we obtain 
\begin{align*}
\mathrm{B} & = \int_{[0,\infty)^3} \ew(X_0X_s)\ew(X_tX_v) \, P_{(T_l,T_m,T_n)} (\rmd (s,t,v)) \\
&=\sigma^4 \int_{[0,\infty)^3} \int_\R f(u)f(u+s)\rmd u \int_\R f(w+t)f(w+v) \rmd w \, P_{(T_l,T_m,T_n)} (\rmd (s,t,v)) \\
&= \sigma^4 \ \ew \bigg[ \int_\R f(u) f(u+T_{l}) \rmd u\int_\R f(u+T_m) f(u+T_{n}) \rmd u \bigg] \, .
\end{align*}
Likewise
\begin{align*}
\mathrm{C} 
&=\sigma^4 \ \ew\bigg[ \int_\R f(u)f(u+T_{m})\rmd u \int_\R f(w+T_{l})f(w+T_{n}) \rmd w \bigg] \, , \quad \text{and} \\
\mathrm{D} &=\sigma^4 \ \ew\bigg[ \int_\R f(u)f(u+T_{n})\rmd u \int_\R f(w+T_{l})f(w+T_{m}) \rmd w \bigg] \, .
\end{align*}
With the definition of $F(s,t)$, the assertion follows.

(b) Observe that, since $P_{\sum_{i=m+1}^{n} W_i} = P_{T_{n-m}}$, by independence of the sequence $W$,
\clearpage
\begin{align*}
& \ew [ F(0,T_l)F(T_m,T_n)] \\
&\quad= \int_{[0,\infty)^3} \int_\R f(u)f(u+s) \rmd u\int_\R f(w+s+t)f(w+s+t+v) \rmd w \\
&\hspace{7cm} P_{\sum_{i=m+1}^{n} W_i} (\rmdd v) P_{\sum_{i=l+1}^{m} W_i} (\rmdd t) P_{T_{l}} (\rmdd s) \\
&\quad= \int_{[0,\infty)} \int_\R f(u)f(u+s) \rmd u \, P_{T_{l}} (\rmdd s) \int_{[0,\infty)} \int_\R f(w)f(w+v) \rmd w \, P_{T_{n-m}} (\rmdd v)  \\ 
&\quad = \ew(F(0,T_l))\ew(F(0,T_{n-m})) \, ,
\end{align*}
which gives the result.
\qed

\vspace*{0.2cm}

From Lemma \ref{lem4.5}, the following proposition gives the expression of $n\cov(\gamma^*_n(p),\gamma_n^*(q))$ as $n \to \infty$ needed in the upcoming central limit theorem.

\begin{proposition}\label{prop4.7}
Let $L=(L_t)_{t \in \R}$ be a L\'evy process with expectation zero and finite fourth moment, and denote $\sigma^2:=\ew(L_1^2)$ and $\eta:=\sigma^{-4}\ew(L_1^4)$. Suppose that $f \in L^2(\R) \cap L^4(\R)$, and let $X$ and $Y$ be defined by (\ref{eq4.1}) and (\ref{eq1.3}) with $(T_n)_{n \in \Z}$ by (\ref{eq1.2}). 
Denote
\begin{align*}
F(s,t) &:= \int_\R f(u+s)f(u+t) \rmd u \, , \quad s,t \in \R \, , \quad \text{and} \\
\kappa_f(k,l,m) &:= (\eta-3) \sigma^4\int_\R f(u) \ew(f(u+T_{k})f(u+T_{l})f(u+T_{m})) \rmd u  \\
&\ + \sigma^4 \ew(F(0,T_l)F(T_k,T_m)) + \sigma^4 \ew(F(0,T_m)F(T_k,T_l)) \, , \quad \text{for } k,l,m \in \Z \, .
\end{align*}
Let $p,q \in \N_0$, denote $Z_{p,i}:=Y_iY_{i+p}$, $Z_{q,j}:=Y_jY_{j+q}$ for $i,j \in \Z$ and assume that 
\begin{align}\label{eq4.9}
\int_\R |f(u)| \sum_{k\in\Z} \ew|f(u+T_p)f(u+T_{k})f(u+T_{k+q})| \rmd u < \infty \, ,
\end{align}
and
\begin{align}\label{eq4.10}
\sum_{k\in\Z} \ew\bigg[ \bigg( \int_\R |f(u) f(u+T_{k})|\rmd u \bigg)^2 \bigg] < \infty \, .
\end{align}
Then
\begin{align}
&\cov(Z_{p,i},Z_{q,j}) = \kappa_f(p,j-i,j-i+q) + \sigma^4 \cov(F(0,T_p),F(T_{j-i},T_{j-i+q})) \, , \label{eq4.11} \\
&\cov(F(0,T_p),F(T_{j-i},T_{j-i+q})) = 0 \quad \text{for } j-i \leq p \ \text{ or } \ j-i \leq q \, , \label{eq4.12} \\
&\sum_{k \in \Z} |\cov(Z_{p,0},Z_{q,k})| < \infty \, , \quad \sum_{k\in\Z} |\kappa_f(p,k,k+q)| < \infty \, , \label{eq4.13} 
\end{align}
and
\begin{align}\label{eq4.14}
\begin{split}
\lim_{n \to \infty} n \, \cov(\gamma_n^*(p),\gamma_n^*(q)) &= \sum_{k \in \Z} \cov(Z_{p,0},Z_{q,k}) = \\
&\hspace{-1cm}\sum_{k\in \Z} \kappa_f(p,k,k+q) + \sigma^4 \sum_{k=-q+1}^{p-1} \cov(F(0,T_{p}),F(T_{k},T_{k+q}))  \, .
\end{split}
\end{align}
where $\gamma_n^*(p)$ and $\gamma_n^*(q)$ are defined in (\ref{eq4.2}).
\end{proposition}

\proof From Lemma \ref{lem4.5} (a), since $\ew(Y_iY_{i+p}) = \sigma^4 \ew(F(T_i,T_{i+p})) = \sigma^4 \ew(F(0,T_p))$, and by the stationarity of $Y$, we have
\begin{align*}
\cov(Z_{p,i},Z_{q,j}) &= \ew(Y_iY_{i+p}Y_jY_{j+q}) - \ew(Y_iY_{i+p})\ew(Y_jY_{j+q}) \\
&=(\eta-3)\sigma^4 \ew\bigg(\int_\R f(u) f(u+T_p)f(u+T_{j-i})f(u+T_{j-i+q}) \rmd u \bigg) \nn \\
&\quad + \sigma^4 \ew(F(0,T_p)F(T_{j-i},T_{j-i+q})) +\sigma^4 \ew(F(0,T_{j-i})F(T_p,T_{j-i+q})) \nn \\
&\quad + \sigma^4 \ew(F(0,T_{j-i+q})F(T_p,T_{j-i})) - \sigma^4\ew(F(0,T_p))\ew(F(T_{j-i},T_{j-i+q})) \nn \\
&= \kappa_f(p,j-i,j-i+q) + \sigma^4 \cov(F(0,T_p),F(T_{j-i},F_{j-i+q}))
\end{align*}
which is (\ref{eq4.11}). Equation (\ref{eq4.12}) is an immediate consequence of Lemma \ref{lem4.5} (b). For the proof of (\ref{eq4.13}), by (\ref{eq4.11}) and (\ref{eq4.12}), it is enough to show that $\sum_{k\in\Z} |\kappa_f(p,k,k+q)| < \infty$. To see this, observe that
\begin{align*}
\sum_{k\in\Z} |\kappa_f(p,k,k+q)| &\leq (\eta-3) \sigma^4 \int_\R |f(u)| \sum_{k \in \Z} |\ew(f(u+T_p)f(u+T_k)f(u+T_{k+q}))| \rmd u \\
&\quad + \sigma^4 \sum_{k\in\Z} |\ew(F(0,T_{k})F(T_p,T_{k+q}))| + \sigma^4 \sum_{k\in\Z} |\ew(F(0,T_{k+q})F(T_p,T_{k}))| \, .
\end{align*}
The first of these summands is finite by (\ref{eq4.9}) and the second is finite since, by using Cauchy-Schwarz's inequality twice, 
\begin{align*}
\sum_{k\in\Z} |\ew(F(0,T_{k})F(T_p,T_{k+q}))| 
\leq \bigg(\sum_{k\in\Z} \ew(F(0,T_k))^2\bigg)^{1/2} \bigg(\sum_{k\in\Z} \ew(F(T_p,T_{k+q}))^2\bigg)^{1/2}
\end{align*}
which is finite by (\ref{eq4.10}). The same argument yields finiteness of the third summand, showing (\ref{eq4.13}). 

To see (\ref{eq4.14}), observe that by the stationarity of $Y$, with $k=j-i$,
\begin{align*}
n\cov(\gamma_n^*(p),\gamma_n^*(q)) 
= \sum_{i,j=1}^n \frac{1}{n} \cov(Z_{p,0},Z_{q,j-i}) 
= \sum_{k=-n+1}^{n-1} \frac{n-|k|}{n} \cov(Z_{p,0},Z_{q,k}) \, .
\end{align*}
Since $\sum_{k \in \Z} |\cov(Z_{p,0},Z_{q,k})| < \infty$ by (\ref{eq4.13}), the latter converges to $\sum_{k \in \Z} \cov(Z_{p,0},Z_{q,k})$ as $n\to \infty$ by the dominated convergence theorem, which together with (\ref{eq4.11}) and (\ref{eq4.12}) finishes the proof of (\ref{eq4.14}).
\qed

\begin{remark}
(a) If $q=0$ or $p=0$, it is easy to see that $\cov(F(0,T_p)F(T_k,T_{k+q}))=0$ for all $k \in \{-q+1,\dots,p-1\}$. Hence, the second summand in (\ref{eq4.14}) disappears. 
 
\noindent(b) If we choose $(T_n)_{n \in \Z}$ to be deterministic, i.e. $T_n = n\Delta$ for $n \in \Z$ and some $\Delta > 0$, it is easy to see that (\ref{eq4.9}) and (\ref{eq4.10}) are implied by assumptions (3.2) and (3.11) of Cohen and Lindner \cite{CL:13} to establish the asymptotic normality of the sample autocovariance of the moving average process sampled on a lattice. (\ref{eq4.9}) then reduces to (3.10) of \cite{CL:13}, which was shown to be implied by (3.2) of \cite{CL:13}.
\end{remark}

\begin{remark}\label{rem4.9}
Similarly to Remark \ref{rem3.5}, a sufficient condition for the validity of (\ref{eq4.9}) and (\ref{eq4.10}) is that $|f(u)|\leq K (|u|^{-\alpha}\land 1)$ for some $K>0$ and $\alpha > 1/2$.

To see this, observe that, by (\ref{eq3.14}),
$$\sum_{k \in \Z} \ew \bigg( \int_\R |f(u)f(u+T_k)|\bigg)^2 \leq \sum_{k\in\Z} C_\alpha^2 \ew(|T_k|^{-2\alpha}\land 1) < \infty \, ,$$
by the the same calculations as in Remark \ref{rem3.5}. Hence, (\ref{eq4.10}) is true. To establish (\ref{eq4.9}), observe that
for some $C_{2\alpha}$, by (\ref{eq3.14}) and applying the Cauchy-Schwarz inequality twice, 
\begin{align*}
&\sum_{k\in\Z} \ew \bigg(\int_\R |f(u)f(u+T_{k})| |f(u+T_p)f(u+T_{k+q})| \rmd u \bigg) \\
&\quad \leq \sum_{k\in\Z} C_{2\alpha} (\ew(|T_k|^{-2\alpha} \land 1))^{1/2} (\ew(|T_{k+q}-T_p|^{-2\alpha} \land 1))^{1/2} \\
&\quad \leq  C_{2\alpha} \bigg(\sum_{k\in\Z} \ew(|T_k|^{-2\alpha} \land 1) \bigg)^{1/2} \bigg( \sum_{k\in\Z} \ew(|T_{k+q}-T_p|^{-2\alpha} \land 1)\bigg)^{1/2} \\
&\quad = C_{2\alpha} \sum_{k\in\Z} \ew(|T_k|^{-2\alpha} \land 1) \, ,
\end{align*}
and the latter is finite by the calculation in Remark \ref{rem3.5}.
\end{remark}

The next proposition shows that similar results as obtained in Proposition \ref{prop4.7} are valid for the truncated sequence $Y^{(m)}$.

\begin{proposition}\label{prop4.10}
Let the assumptions and notations of Proposition \ref{prop4.7} be satisfied. For $m \in \N$, define $f_m:=f\ind{[-m/2,m/2]}$, $F_m:=\int_\R f_m(u+s)f_m(u+t)$ for $s,t \in \R$, $X_t^{(m)}:=\int_\R f_m(t-u)\rmd L_u$, $Y_n^{(m)}:=X_{T_n}^{(m)}$. Let $p,q \in \N_0$ and 
$$\gamma^{*,m}_n(h) :=\frac{1}{n} \sum_{k=1}^n Y_k^{(m)}Y_{k+h}^{(m)} \, , \quad h=p,q \, .$$
Then (\ref{eq4.9}) and (\ref{eq4.10}) also hold for $f_m$, and for all $k \in \Z$
\begin{align}
|\kappa_{f_m}(p,k,k+q)| &\leq \kappa_{|f|}(p,k,k+q) \, , \nn \\
\kappa_{f_m}(p,k,k+q) &\to \kappa_{f}(p,k,k+q) \quad  \text{as } \ m \to \infty \, , \label{eq4.15} \\
\cov(F_m(0,T_p)F_m(T_k,T_{k+q})) &\to \cov(F(0,T_p)F(T_k,T_{k+q})) \quad \text{as } \ m \to \infty \, , \label{eq4.16} 
\end{align}
and
\begin{align}\label{eq4.17}
\begin{split}
&\lim_{m\to\infty}\lim_{n\to\infty} n \cov(\gamma^{*,m}_n(p),\gamma^{*,m}_n(q)) \\
&\hspace{2cm}= \sum_{k\in \Z} \kappa_f(p,k,k+q) + \sigma^4 \sum_{k=-q+1}^{p-1} \cov(F(0,T_{p}),F(T_{k},T_{k+q}))  \, .
\end{split}
\end{align}
\end{proposition}

\proof That (\ref{eq4.9}) and (\ref{eq4.10}) also hold for $f_m$ is clear since $|f_m| \leq |f|$, as is $|\kappa_{f_m}| \leq \kappa_{|f|}$. Since $|f_m|\leq |f|$ and $f_m \to f$ as $m \to \infty$, the dominated convergence theorem shows (\ref{eq4.15}) and (\ref{eq4.16}). And (\ref{eq4.17}) then follows from (\ref{eq4.15}), (\ref{eq4.16}), and (\ref{eq4.14}) again by the dominated convergence theorem.
\qed

\vspace{0.2cm}

Now, we can establish the multivariate asymptotic normality of the sample autocovariance and sample autocorrelation.

\begin{theorem}\label{thm4.12}
Let $L=(L_t)_{t \in \R}$ be a L\'evy process with expectation zero and such that $\ew(|L_1|^4 (\log^+|L_1|)^2) < \infty$. Denote $\sigma^2:=\ew(L_1^2)$ and $\eta:=\sigma^{-4}\ew(L_1^4)$. Let $h \in \N_0$, suppose that $f \in L^2(\R) \cap L^{4}(\R)$, $\int_\R |f(s)|^4 (\log^+|f(s)|)^2 \rmd s < \infty$, and assume that (\ref{eq4.9}) and (\ref{eq4.10}) hold for all $p,q \in \{0,\dots,h\}$. 

\noindent(a) Then
\begin{align*}
\sqrt{n} (\gamma_n^*(0)-\gamma(0),\dots,\gamma_n^*(h)-\gamma(h))' \dcon N(0,\mathbf{Z}) \, , \quad n \to \infty \, ,
\end{align*}
where $\gamma(h)=\ew(Y_0Y_h)$ and $\mathbf{Z}=(\mathbf{Z}_{pq})_{p,q=0,\dots,h} \in \R^{(h+1) \times (h+1)}$ is the covariance matrix defined by
\begin{align*}
\mathbf{Z}_{pq} &= \sigma^4 \sum_{k=-q+1}^{p-1} \cov(F(0,T_{p}),F(T_{k},T_{k+q})) + \sum_{k \in \Z} \kappa_f(p,k,k+q)
\end{align*}
with $\kappa_f(p,k,k+q)$ for $k,p,q \in \Z$ and $F(s,t)$ given as in Proposition \ref{prop4.7}. 
\vspace*{0cm}

\noindent (b) If additionally 
\begin{align}\label{eq4.19a}
\int_\R |f(u)| \sum_{k\in\Z} \ew|f(T_k+u)| \rmd u < \infty
\end{align}
hold, and we denote by
$$\widehat{\gamma}_n(j) =\frac{1}{n} \sum_{k=1}^{n-j} (Y_k - \overline{Y}_n)(Y_{k+j} - \overline{Y}_n) \, , \quad j=0,1,\dots,n-1 \, ,$$
the sample autocovariance, then we have for each $h \in \N_0$
\begin{align*}
\sqrt{n} (\widehat{\gamma}_n(0)-\gamma(0),\dots,\widehat{\gamma}_n(h)-\gamma(h))' \dcon N(0,\mathbf{Z}) \, , \quad n \to \infty \, ,
\end{align*}
where $\mathbf{Z}$ as defined in (a).

\noindent (c) Let $\rho_n^*(p)=\gamma^*_n(p)/\gamma_n^*(0)$ and $\widehat{\rho}_n(p)=\widehat{\gamma}_n(p)/\widehat{\gamma}_n(0)$ for $p \in \N$. Suppose that $f \neq 0$ on a set of positive Lebesgue measure. Then, under the assumptions of (a), we have for each $h \in \N$
\begin{align*}
\sqrt{n} (\rho_n^*(1)-\rho(1),\dots,\rho_n^*(h)-\rho(h))' \dcon N(0,\mathbf{W}) \, , \quad n \to \infty \, ,
\end{align*}
where $\mathbf{W}=(\mathbf{W}_{pq})_{p,q=1,\dots,h} \in \R^{h \times h}$ is given by
$$\mathbf{W}_{pq} =(\mathbf{Z}_{pq}-\rho(p)\mathbf{Z}_{0q}-\rho(q)\mathbf{Z}_{p0}+\rho(p)\rho(q)\mathbf{Z}_{00})/\gamma(0)^2 \, .$$
If additionally (\ref{eq4.19a}) is satisfied, then it also holds
\begin{align*}
\sqrt{n} (\widehat{\rho}_n(1)-\rho(1),\dots,\widehat{\rho}_n(h)-\rho(h))' \dcon N(0,\mathbf{W}) \, , \quad n \to \infty \, .
\end{align*}
\end{theorem}

\proof (a) Let $Y^{(m)}$ be as in (\ref{eq2.2}) and (\ref{eq2.1}) with $\mu=0$ and $Z^{(m)}_{p,k}=Y^{(m)}_kY^{(m)}_{k+p}$ for $k \in \Z$ and $p \in \{0,\dots,h\}$. Define
$Q_k := (Z^{(m)}_{0,k},Z^{(m)}_{1,k},\dots,Z^{(m)}_{h,k})' \in \R^{h+1}$.
Then $(Q_k)_{k \in \Z}$ is obviously strictly stationary and we have 
$$\frac{1}{n} \sum_{k=1}^n Q_k = (\gamma^{*,m}_n(0),\dots,\gamma^{*,m}_n(h))' \, ,$$
where $\gamma^{*,m}_n(h)$ as in Proposition \ref{prop4.10}.

By the assumptions on $L$ and $f$, we obtain, by Lemma \ref{lem2.2} (c), $\ew(|Y^{(m)}_0|^{4}(\log^+|Y_0^{(m)}|)^2)=\ew(|X^{(m)}_0|^{4}(\log^+|X_0^{(m)}|)^2)<\infty$, and so $\ew(|Z_{h,0}|^{2} \log^+|Z_{h,0}|)<\infty$, by the Cauchy-Schwarz inequality. Therefore also $\ew(|\lbd'Q_0|^{2}\log^+|\lbd'Q_0|)<\infty$ for all $\lbd \in \R^{h+1}$.

Observe that $(\lbd'Q_n)_{n \in \Z}$ is strongly mixing for each $\lbd \in \R^{h+1}$ with $\alpha_n^{\lbd'Q}\leq \alpha_{n-h}^{Y^{(m)}}$ for all $n > h$, by Remark 1.8 (b) of Bradley \cite{brad:07-1}, such that $(\alpha_n^{\lbd'Q})$ is exponentially decreasing.  
Hence, by Corollary 10.20 (c) of Bradley \cite{brad:07-1} it follows
$$\sqrt{n}\bigg( \frac{1}{n} \sum_{k=1}^n \lbd'Q_k - \lbd' (\gamma^{m}_n(0),\dots,\gamma^{m}_n(h))' \bigg) \dcon N(0,\lbd'\mathbf{Z}^m \lbd) \quad \forall \, \lbd \in \R^{h+1} \, .$$
By the Cram\'er-Wold theorem, we deduce that 
$$\sqrt{n} (\gamma_n^{*,m}(0)-\gamma^m(0),\dots,\gamma_n^{*,m}(h)-\gamma^m(h))' \dcon \mathbf{V}_m \quad \text{as } n \to \infty \, .$$
Here $\mathbf{V}_m \overset{d}{=} N(0,\mathbf{Z}^m)$, where $\mathbf{Z}^m=(\mathbf{Z}^m_{pq})_{p,q=0,\dots,h} \in \R^{(h+1) \times (h+1)}$ is given by 
$$\mathbf{Z}^m_{p,q} = \sum_{k\in \Z} \kappa_{f_m}(p,k,k+q) + \sigma^4 \sum_{k=-q+1}^{p-1} \cov(F_m(0,T_{p}),F_m(T_{k},T_{k+q})) \, ,$$
with $\kappa_{f_m}$ and $F_m$ as in Proposition \ref{prop4.10}.

Also by Proposition \ref{prop4.10}, $\lim_{m \to \infty} \mathbf{Z}^m = \mathbf{Z}$, where Proposition \ref{prop4.7} gives the form and finiteness of $\mathbf{Z}_{pq}$, the entries of $\mathbf{Z}$. Henceforth, $\mathbf{V}_m \dcon \mathbf{V}$ as $m \to \infty $,
where $V \overset{d}{=} N(0,\mathbf{Z})$. 

By Proposition 6.3.9 of Brockwell and Davis \cite{brock:06}, the claim will follow if we can show that
$$\lim_{m \to \infty} \limsup_{n \to \infty} P(n^{1/2}|\gamma^{*,m}_n(p)-\gamma^m(p)-\gamma^*_n(p) + \gamma(p)| > \varepsilon) = 0 \quad \forall \, \varepsilon > 0 \, , \quad p \in \{0,\dots,h\} \, .$$
Since $\ew(\gamma_n^{*,m}(p))=\gamma^{m}(p)$ and $\ew(\gamma_n^{*}(p))=\gamma(p)$, this will follow from Chebychef's inequality if we can show that 
\begin{align*}
\lim_{m\to \infty} \lim_{n \to \infty} \var(n^{1/2} (\gamma^{*}_n(p) -\gamma^{*,m}_n(p)) &= \lim_{m \to \infty} \lim_{n \to \infty} \big[n \var(\gamma^{*}_n(p)) + n \var(\gamma^{*,m}_n(p)) \\
&\hspace{0.5cm}- 2n \cov(\gamma^{*}_n(p),\gamma^{*,m}_n(p)) \big] =0 \quad \forall \, p \in \{0,\dots,h\} \, .
\end{align*}
But since 
$$\lim_{m \to \infty} \lim_{n \to \infty} n\var(\gamma_n^{*,m}(p)) =\lim_{n \to \infty} n \var(\gamma_n^{*}(p)) = \mathbf{Z}_{pp} \, ,$$
by Proposition \ref{prop4.10}, it remains only to show that 
\begin{align}\label{eq4.19}
\lim_{m \to \infty} \lim_{n \to \infty} n \cov(\gamma^*_n(p),\gamma^{*,m}_n(p)) = \mathbf{Z}_{pp} \quad \forall \, p \in \{0,\dots,h\} \, .
\end{align}

In doing so, denote $G_m(s,t):=\int_\R f(u+s)f_m(u+t) \rmd u$ and $F_m(s,t):=\int_{\R} f_m(u+s)f_m(u+t) \rmd u$. Observe first that from Lemma \ref{lem4.3} (b), similar to the proof of Lemma \ref{lem4.5} (a), by conditioning on $T_p$, $T_k$, and $T_{k+p}$, that for $k \in \Z$, we have
\clearpage
\begin{align*}
\cov(Z_{p,0},Z_{p,k}^{(m)}) &= \ew(Y_0Y_{p}Y_k^{(m)}Y_{k+p}^{(m)}) - \ew(Y_0Y_p)\ew(Y_k^{(m)}Y_{k+p}^{(m)}) \\
&=(\eta-3)\sigma^4 \ew\bigg(\int_\R f(u) f(u+T_p)f_m(u+T_{k})f_m(u+T_{k+p})  \rmd u \bigg)  \\
&\quad + \sigma^4 \ew(F(0,T_p)F_m(T_{k},T_{k+p})) +\sigma^4 \ew(G_m(0,T_{k})G_m(T_p,T_{k+p}))  \\
&\quad + \sigma^4 \ew(G_m(0,T_{k+p})G_m(T_p,T_{k})) - \sigma^4\ew(F(0,T_p))\ew(F_m(T_{k},T_{k+p})) \, . 
\end{align*}
Further, as in the proof of Lemma \ref{lem4.5} (b), it follows that
$$\ew(F(0,T_p)F_m(T_{k},T_{k+p})) = \ew(F(0,T_p))\ew(F_m(T_{k},T_{k+p})) \quad \text{ when } |k| \geq p \, .$$
Denoting
\begin{align*}
\kappa_{f,f_m}(p,k,k+p) &:= (\eta-3) \sigma^4 \ew\bigg(\int_\R f(u) f(u+T_{p})f_m(u+T_{k})f_m(u+T_{k+p}) \rmd u \bigg)  \\
&\ + \sigma^4 \ew(G_m(0,T_k)G_m(T_p,T_{k+p})) + \sigma^4 \ew(G_m(0,T_{k+p})G_m(T_p,T_k)) \, ,
\end{align*}
\noindent we hence have
\begin{align*}
\cov(Z_{p,0},Z_{p,k}^{(m)}) = \begin{cases}
\kappa_{f,f_m}(p,k,k+p) \, , &\quad |k| \geq p \, , \\
\kappa_{f,f_m}(p,k,k+p) + \sigma^4 \cov(F(0,T_p)F_m(T_k,T_{k+p})) \, , &\quad |k| < p \, .
\end{cases}
\end{align*}
Next, observe that as in the proof of Proposition \ref{prop4.10}, since $|f_m|\leq |f|$, for all $k \in \Z$,
\begin{align*}
|\kappa_{f,f_m}(p,k,k+p)| &\leq \kappa_{|f|}(p,k,k+p) \quad \forall \, m \in \N \, ,  \\
\kappa_{f,f_m}(p,k,k+p) &\to \kappa_{f}(p,k,k+p) \quad  \text{as } \ m \to \infty \, ,  \\
\cov(F(0,T_p)F_m(T_k,T_{k+p})) &\to \cov(F(0,T_p)F(T_k,T_{k+p})) \quad \text{as } \ m \to \infty \, . 
\end{align*}
By stationarity, we obtain for $n \geq p$
\begin{align*}
&n \cov(\gamma^*_n(p),\gamma^{*,m}_n(p)) = \frac{1}{n} \sum_{i,j=1}^n \cov(Z_{p,0},Z_{p,j-i}^{(m)}) = \sum_{k=-n+1}^{n-1} \frac{n-|k|}{n} \cov(Z_{p,0},Z_{p,k}^{(m)}) \\
&\quad= \sum_{k=-n+1}^{n-1} \frac{n-|k|}{n} \kappa_{f,f_m}(p,k,k+p) + \sigma^4 \sum_{k=-p+1}^{p-1} \frac{n-|p|}{n}\cov(F(0,T_p)F_m(T_k,T_{k+p})) \, .
\end{align*}
Applying Lebesgue's dominated convergence theorem once then gives
\begin{align*}
\lim_{n \to \infty} & n \cov(\gamma^*_n(p),\gamma^{*,m}_n(p)) \\
&\quad= \sum_{k=-\infty}^{\infty} \kappa_{f,f_m}(p,k,k+p) + \sigma^4 \sum_{k=-p+1}^{p-1} \cov(F(0,T_p)F_m(T_k,T_{k+p})) \, ,
\end{align*}
and applying it a second time gives
\begin{align*}
\lim_{m \to \infty} \lim_{n \to \infty} & n \cov(\gamma^*_n(p),\gamma^{*,m}_n(p)) \\
&\quad= \sum_{k=-\infty}^{\infty} \kappa_{f}(p,k,k+p) + \sigma^4 \sum_{k=-p+1}^{p-1} \cov(F(0,T_p)F(T_k,T_{k+p})) \, ,
\end{align*}
which is (\ref{eq4.19}). This finishes the proof of (a).

\clearpage

(b) This follows if we can show that $\sqrt{n}|\gamma_n^{*}(p)-\widehat{\gamma}_n(p)| \to 0$ in probability for $n \to \infty$ and $p \in \{0,\dots,h\}$. The latter can be done in exactly the same way as in the proof of Proposition 7.3.4 in Brockwell and Davis \cite{brock:06} with $X$ replaced by $Y$ in connection with the observation that, by Theorem \ref{thm2.6}, $\sqrt{n}\overline{Y}_n$ converges in distribution to a normal random variable as $n \to \infty$, and hence $\overline{Y}_n$ must converges to $0$ in probability as $n \to \infty$.

(c) Follows readily as in the proof of Theorem 7.2.1 in Brockwell and Davis \cite{brock:06}.
\qed

\begin{remark}
(a) Due to the form of $\mathbf{Z}$, there seems to be no simplification for $\mathbf{W}$ possible. Also observe that $\mathbf{W}$ in general depends on $\eta$ as seen in Theorem 3.5 (c) of Cohen and Lindner \cite{CL:13}. \\
\noindent(b) Part (a) of Theorem \ref{thm4.12} in particular applies if $|f(u)|\leq K (|u|^{-\alpha}\land 1)$ for some $K>0$ and $\alpha > 1/2$ which can be seen by Remark \ref{rem4.9}. \\
\noindent(c) Similarly, part (b) of Theorem \ref{thm4.12} applies if $|f(u)|\leq K (|u|^{-\alpha}\land 1)$ for some $K>0$ and $\alpha > 1$ as shown in Remark \ref{rem3.5}.
\end{remark}

\section{An Application to Parameter Estimation of the Ornstein-Uhlenbeck Process}
\setcounter{equation}{0}

In this section, we present a parameter estimation of a L\'evy driven Ornstein-Uhlenbeck (OU) process sampled at a Poisson process. An OU process is a continuous time moving average process $X=(X_t)_{t \in \R}$ with kernel function $f \colon \R \to \R , s \mapsto \e^{-as}\ind{[0,\infty)}(s)$ and mean reverting parameter $a > 0$. This yields $X_t = \int_{-\infty}^t \e^{-a(t-s)} \rmd L_s$, where $L=(L_t)_{t \in \R}$ is a L\'evy process with zero mean and $\sigma^2=\ew(L_1^2)<\infty$. We define $Y_n := X_{T_n}$, $n \in \Z$, where $(T_n)_{n \in \Z}$ is given by (\ref{eq1.2}) with $W=(W_n)_{n \in \Z \setminus \{0\}}$ a sequence of i.i.d. random variables independent of $L$ and such that $W_1 \sim \mathrm{Exp}(\lbd)$, $\lbd > 0$. 

By Proposition \ref{prop1.1}, $Y=(Y_n)_{n \in \Z}$ is strictly stationary, $\ew(Y_0)=0$, and $\ew(Y_0^2)<\infty$. Then we obtain, by Proposition \ref{prop2.5},
\begin{align*}
\gamma(h) = \ew(Y_0Y_h) &= \sigma^2 \int_\R f(u) \ew(f(T_h+u)) \rmd u 
= \frac{\sigma^2}{2a} \bigg(\frac{\lbd}{a+\lbd}\bigg)^h 
\end{align*}
as the autocovariance function of the process $Y$. Further $\gamma(0) = \sigma^2/2a$, and the autocorrelation function $\rho(h) = \big(\frac{\lbd}{a+\lbd}\big)^h$. In particular, $\rho(1)=\frac{\lbd}{\lbd+a}$ and we can determine the mean reverting parameter $a$ by 
\begin{align}\label{eq51}
a = \bigg(\frac{1}{\rho(1)}-1 \bigg) \, .
\end{align}
We define an estimator $a^*$ for $a$, assuming $\lbd$ being a known parameter of the distribution of $W$, as
$a^*=\lbd\bigg(\frac{1}{\rho^*(1)} - 1 \bigg)$,
where $\rho^*(1)=\gamma^*(1)/\gamma^*(0)$ with $\gamma^*(h) = \frac{1}{n} \sum_{k=1}^n Y_k Y_{k+h}$. We can then give the following theorem.

\clearpage

\begin{theorem}\label{thm51}
Let $L=(L_t)_{t \in \R}$ be a L\'evy process with mean zero, $\sigma^2=\ew(L_1^2)$ and $\eta=\sigma^{-4}\ew(L_1^4)$, and $\ew(|L_1|^4(\log^+|L_1|)^2)<\infty$, $(T_n)_{n \in \Z}$ be defined as in (\ref{eq1.2}) such that $W_1 \sim \mathrm{Exp}(\lbd)$ for some $\lbd > 0$, and $X=(X_t)_{t \in \R}$ an OU process with parameter $a > 0$. 
Then 
\begin{align*}
\sqrt{n}(a^*-a) \dcon N\bigg(0,\frac{(\lbd+a)^4}{\lbd^2} \mathbf{W}_{11} \bigg) \, , \quad n \to \infty \, ,
\end{align*}
where
\begin{align}\label{*r53}
\mathbf{W}_{11} = \bigg( \frac{\lbd}{\lbd+2a} - \frac{\lbd^2}{(\lbd+a)^2} \bigg) ((\eta-3)a+3) + \frac{2a}{\lbd+2a} \, .
\end{align}
\end{theorem}

\proof As usual, $Y_k:=X_{T_k}$, $k \in \Z$, and $Z_{k,h}:=Y_kY_{k+h}$, $k \in \Z$, $\gamma_n^*(h)=\frac{1}{n}\sum_{k=1}^n Z_{k,h}$, $h=0,\dots,n-1$ and hence $\rho_n^*(1)=\gamma^*_n(1)/\gamma_n^*(0)$. Observe that $P(W_1 > 0)>0$, since $W$ is exponentially distributed and it has positive support. Further, $f \in L^2(\R)\cap L^4(\R)$ is obvious as is $\int_\R |f(s)|^4(\log^+|f(s)|)^2\rmd s < \infty$, and, since clearly $f(s)\leq K(|s|^{-\alpha}\land 1)$ for some $K > 0$ and $\alpha > 1/2$, it follows, by Remark \ref{rem4.9}, that (\ref{eq4.9}) and (\ref{eq4.10}) are satisfied. Therefore, by Theorem \ref{thm4.12} (c), we have
\begin{align*}
\sqrt{n}(\rho_n^*(1)-\rho(1)) \dcon N(0,\mathbf{W}_{11} ) \, , \quad n \to \infty \, ,
\end{align*}
where 
$$\mathbf{W}_{11} =(\mathbf{Z}_{11}-2\rho(1)\mathbf{Z}_{01}+\rho(1)^2\mathbf{Z}_{00})/\gamma(0)^2 = \frac{4a^2}{\sigma^4} \bigg(\mathbf{Z}_{11}-2\frac{\lbd}{a+\lbd}\mathbf{Z}_{01}+\bigg(\frac{\lbd}{a+\lbd}\bigg)^2\mathbf{Z}_{00}\bigg)$$
with $\mathbf{Z}_{pq}$ for $p,q \in \{0,1\}$ given as in Theorem \ref{thm4.12} (a).

Thus, under our assumptions on the distribution of $W$, an easy but tedious calculation yields that $\mathbf{W}_{11}$ is given by (\ref{*r53}).

To complete the proof, define $g\colon \R \to \R , x\mapsto\lbd\big(\frac{1}{x}-1\big)$ such that $g(\rho_n^*(1))=a^*$ and the delta-method, cf. Proposition 6.4.3 in Brockwell and Davis \cite{brock:06}, yields 
\begin{align*}
\sqrt{n}(a^*-a) \dcon N(0,g'(\rho(1)) \mathbf{W}_{11} g'(\rho(1))) \, , \quad n \to \infty \, ,
\end{align*}
where $g'(\rho(1))= -\frac{(\lbd+a)^2}{\lbd}$. 
\qed

\vspace{0.2cm}

Next, we consider the case when the parameter $\lbd$ of $W_1 \sim \mathrm{Exp}(\lbd)$ is unknown. Since, in addition to the observations $Y_1,\dots,Y_{n+1}$, we also have the observation times $T_1,\dots,T_{n+1}$, we also observe the waiting times $W_i=T_i-T_{i-1}$, $i=1,\dots,n+1$, and hence can define 
$\widehat{\lbd} := \bigg(\frac{1}{n} \sum_{k=1}^n W_{k+1} \bigg)^{-1}$,
which by the strong law of large numbers is a strongly consistent estimator for $\lbd$, since $\ew(W_1)=\lbd^{-1}$. 

By (\ref{eq51}), this suggests the estimator 
$\widehat{a}=\widehat{\lbd}\bigg(\frac{1}{\rho^*(1)} - 1 \bigg)$.
Since $\rho^*(1)$ and $\widehat{\lbd}$ are consistent estimators, so is $\widehat{a}$. The asymptotic normality of $\widehat{a}$ is given in the following theorem.

\begin{theorem}
Let $L=(L_t)_{t \in \R}$ be a L\'evy process with mean zero, $\sigma^2=\ew(L_1^2)$, $\eta=\sigma^{-4}\ew(L_1^4)$, and $\ew(|L_1|^4(\log^+|L_1|)^2)<\infty$. Assume that $(T_n)_{n \in Z}$ is defined as in (\ref{eq1.2}) such that $W_1 \sim \mathrm{Exp}(\lbd)$ for some $\lbd > 0$, and $X=(X_t)_{t \in \R}$ is a OU process with parameter $a > 0$. 
Then 
\begin{align*}
\sqrt{n}(\widehat{a}-a) \dcon N\bigg(0,\frac{(\lbd+a)^4}{\lbd^2} \mathbf{W}_{11} - a^2 \bigg) \, , \quad n \to \infty \, ,
\end{align*}
where $\mathbf{W}_{11}$ is given by (\ref{*r53}).
\end{theorem}

\proof 
%
%
For $m\in \N$ define $f_m:=f\ind{[-m/2,m/2]}$ and $Y_n^{(m)}:=\int_\R f_m(T_n-s) \rmd L_s$. Then the sequences $(Y_n^2,Y_nY_{n+1},T_{n+1}-T_n)_{n \in \Z}$ and $((Y_n^{(m)})^2,Y_n^{(m)}Y_{n+1}^{(m)},T_{n+1}-T_n)_{n \in \Z}$ are both strictly stationary by Proposition \ref{prop1.1} and the latter is also strongly mixing with exponentially decreasing mixing coefficients by Proposition \ref{prop2.1}. 

Proceeding exactly as in the proof of Theorem \ref{thm4.12}, i.e. establishing first a central limit theorem for the truncated sequences $\gamma_n^{*,m}(h)$ and then letting $m$ tend to infinity, shows that
\begin{align*}
\sqrt{n} \bigg(\bigg(\gamma_n^*(0),\gamma_n^{*}(1), \frac{1}{n} \sum_{k=1}^n W_{k+1}\bigg)-\bigg(\gamma(0),\gamma(1),\frac{1}{\lbd}\bigg)\bigg) \dcon N(0,\Sigma) \, , \quad n \to \infty \, ,
\end{align*}
where  
\begin{align*}
\Sigma =\sum_{k \in \Z} \begin{pmatrix}
\cov(Y_0^2,Y_k^2) & \cov(Y_0^2,Y_kY_{k+1}) & \cov(Y_0^2,T_{k+1}-T_{k}) \\
\cov(Y_0^2,Y_kY_{k+1}) & \cov(Y_0Y_1,Y_kY_{k+1}) & \cov(Y_0Y_1,T_{k+1}-T_{k}) \\
\cov(Y_0^2,T_{k+1}-T_{k}) & \cov(Y_0Y_1,T_{k+1}-T_{k}) & \cov(T_1,T_{k+1}-T_k)
\end{pmatrix} \, .
\end{align*}
An easy but tedious calculation then shows that 
\begin{align*}
\Sigma = \begin{pmatrix}
\mathbf{Z}_{00} & \mathbf{Z}_{01} & 0 \\
\mathbf{Z}_{10} & \mathbf{Z}_{11} & -\frac{\sigma^2}{2(\lbd+a)^2} \\
0 & -\frac{\sigma^2}{2(\lbd+a)^2} & \frac{1}{\lbd^2}
\end{pmatrix} \, 
\end{align*}
with $\mathbf{Z}_{00}$, $\mathbf{Z}_{01}$, and $\mathbf{Z}_{11}$ as in Theorem \ref{thm4.12} (a).

To complete, we define $g:\R^3 \to \R, (x_1,x_2,x_3) \mapsto \frac{1}{x_3}(\frac{x_1}{x_2}-1)$ such that 
$$\widehat{a} = g\bigg(\gamma_n^*(0),\gamma_n^*(1),\frac{1}{n}\sum_{k=1}^n W_{k+1}\bigg) = \widehat{\lbd} \bigg(\frac{1}{\rho_n^*(1)} - 1\bigg) \, .$$
Henceforth, by the delta-method, cf. Proposition 6.4.3 in Brockwell and Davis \cite{brock:06}, we obtain 
\begin{align*}
\sqrt{n}(\widehat{a}-a) \dcon N(0,(\nabla g(\mu)) \Sigma (\nabla g(\mu))') \, , \quad n \to \infty \, ,
\end{align*}
where 
by a straightforward calculation,
$$(\nabla g(\mu)) \Sigma (\nabla g(\mu))' = \frac{(\lbd+a)^4}{\lbd^2} \mathbf{W}_{11} - a^2 $$
and the result follows.
\qed

\begin{remark}
Note that the shrinking phenomenon observed in the asymptotic variance of the estimator $\widehat{a}$ with respect to the asymptotic variance of $a^*$ depends on the non zero asymptotic covariance between the sample autocovariance $\gamma^*_n(1)$ and the estimator $\widehat{\lbd}$.
\end{remark}

\begin{figure}[!ht]
\centering
\subfigure[Expected waiting time $\Delta=20$.]{\includegraphics[width=7.7cm]{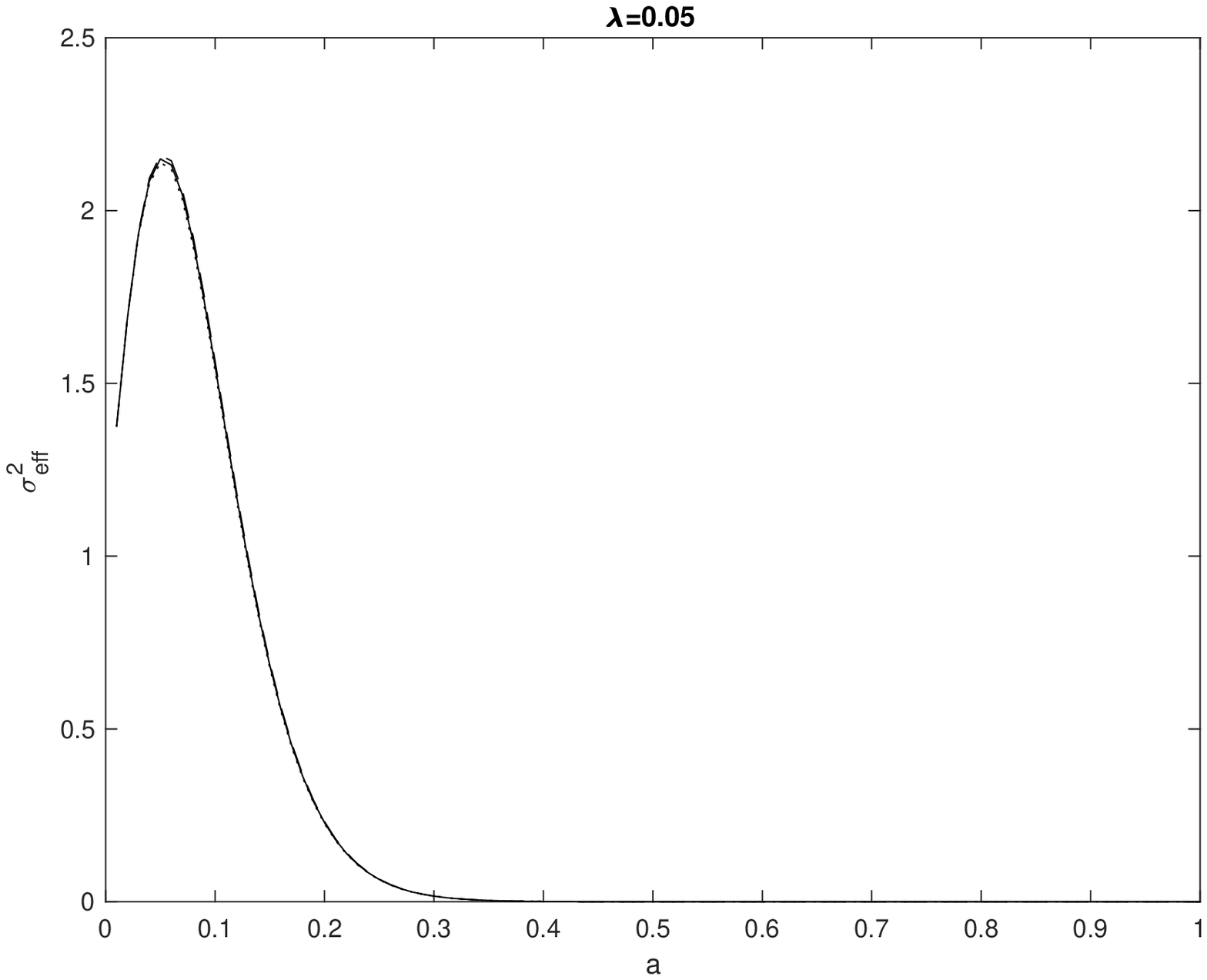}} \hspace{0.33cm}
\subfigure[Expected waiting time $\Delta=2$.]{\includegraphics[width=7.7cm]{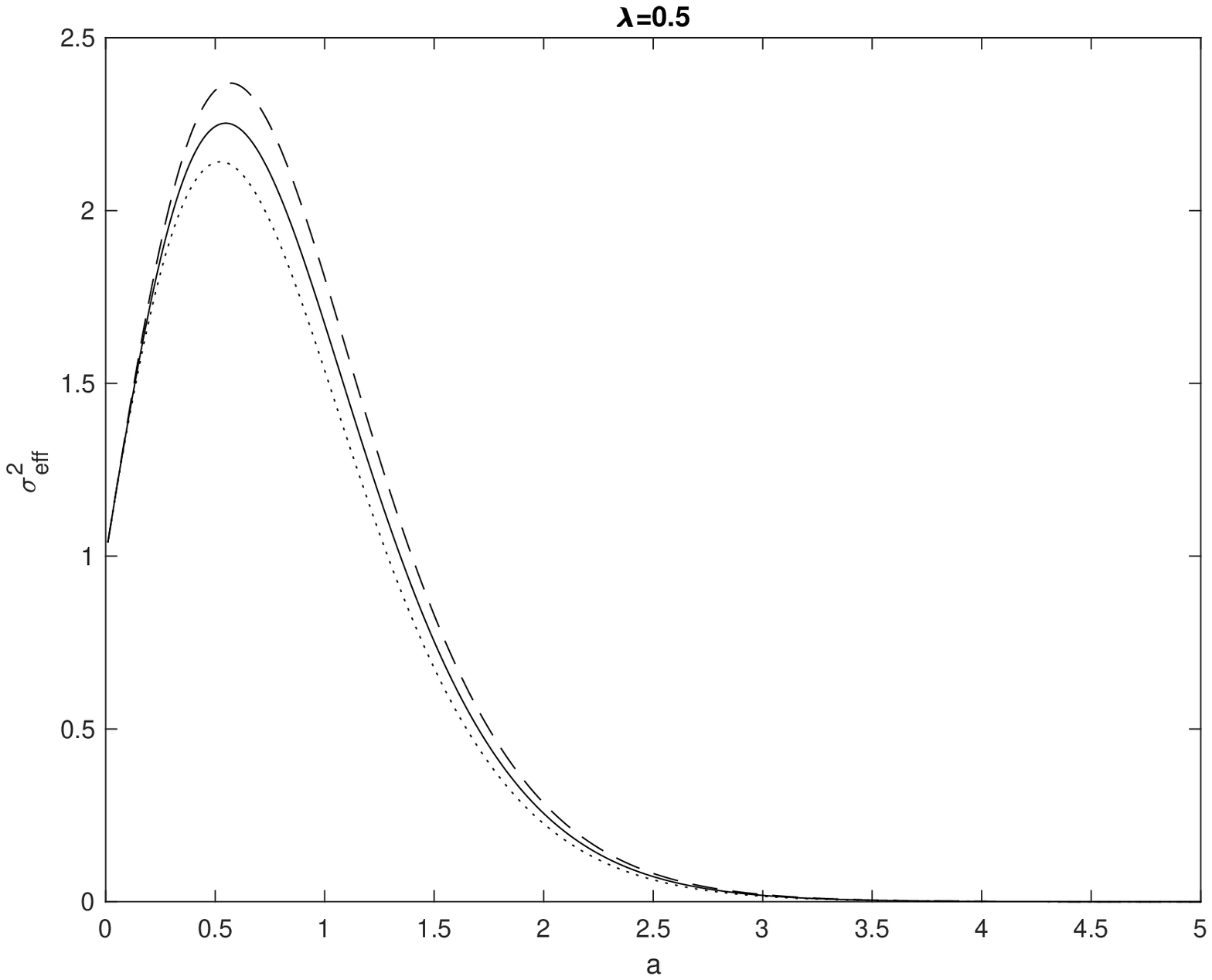}} \vspace{0.35cm}
\subfigure[Expected waiting time $\Delta=1$.]{\includegraphics[width=7.7cm]{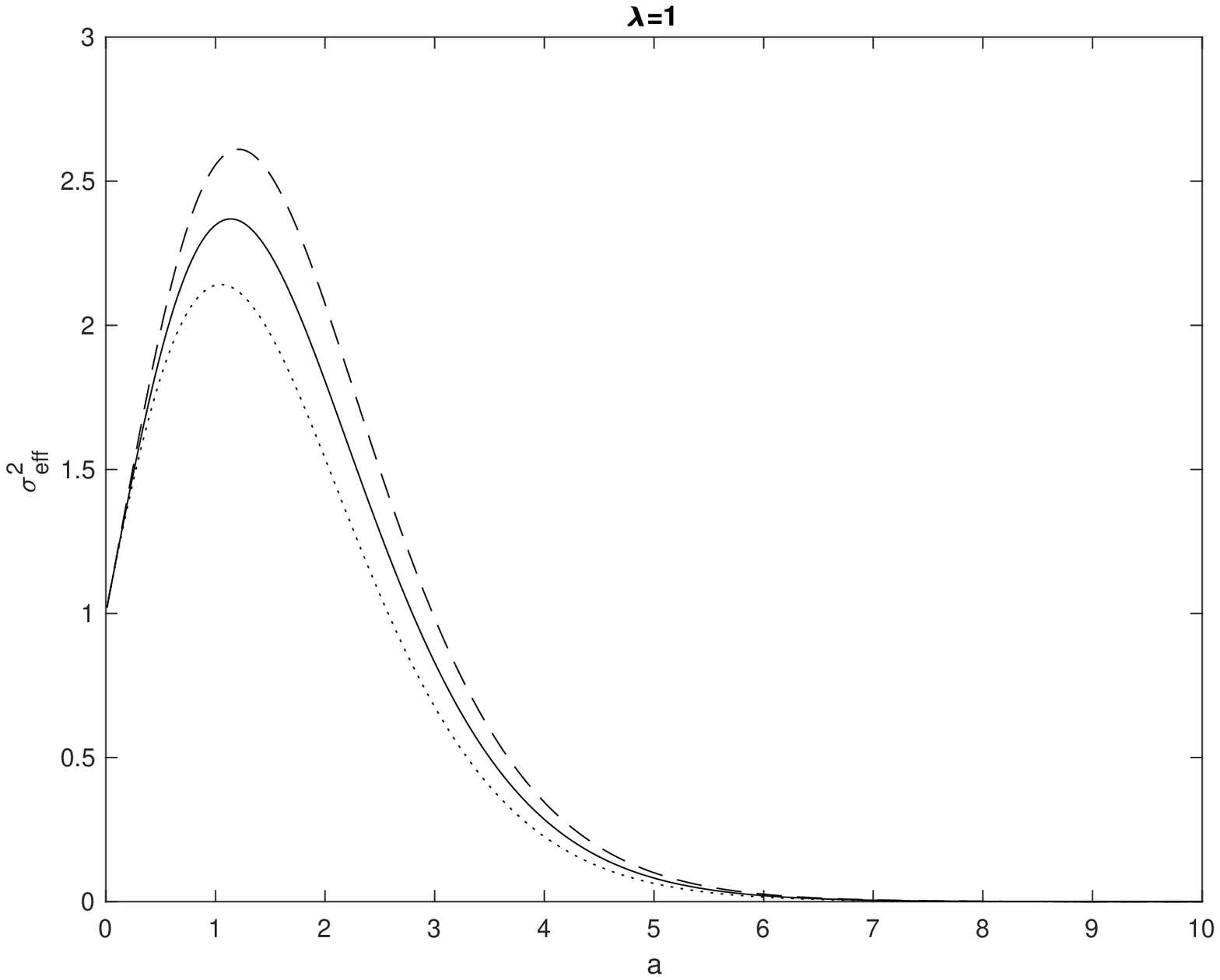}} \hspace{0.33cm}
\subfigure[Expected waiting time $\Delta=0.5$.]{\includegraphics[width=7.7cm]{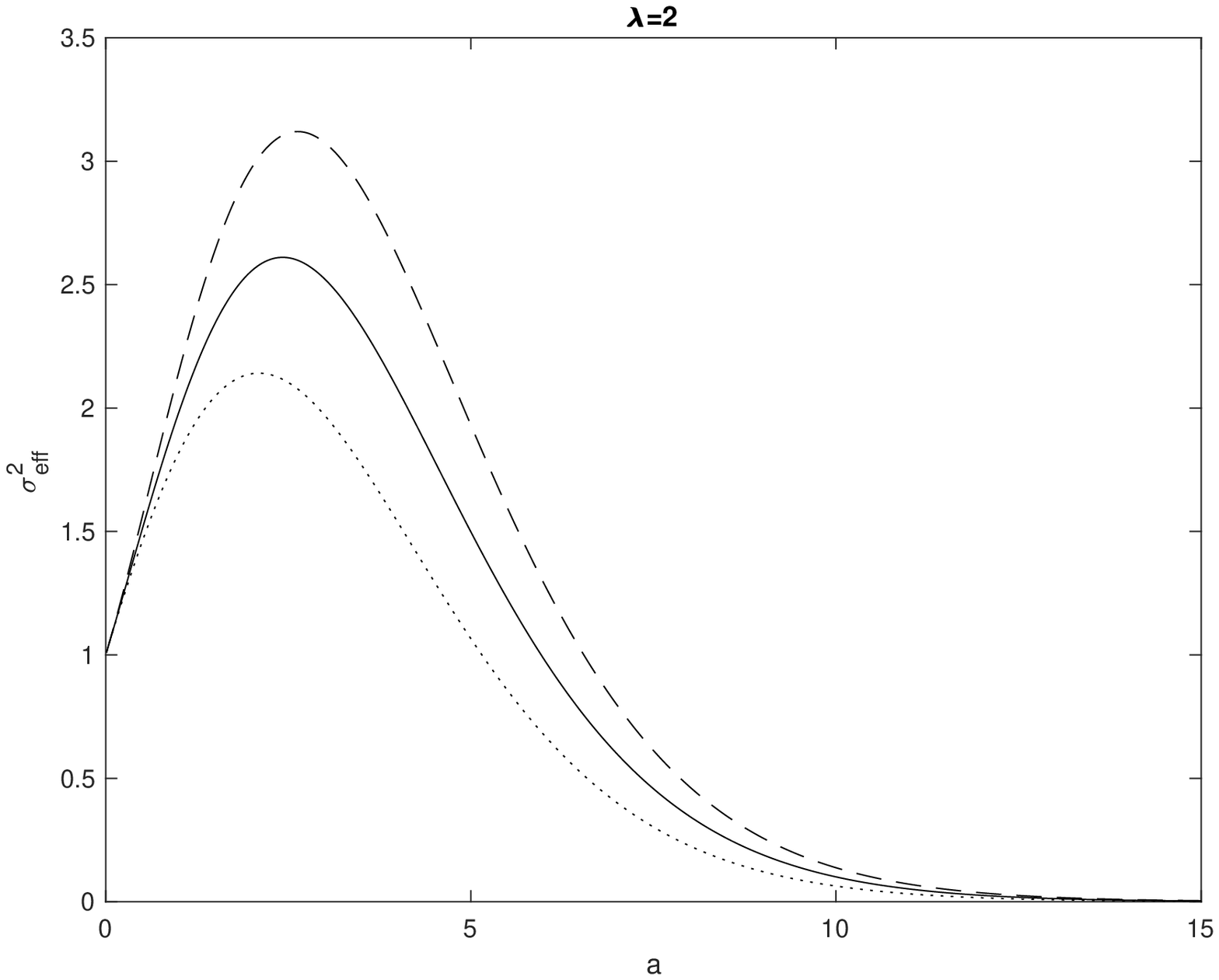}} \hspace{0.33cm}
\caption{$\sigma^2_{eff}$ depending on $a$ and $\lbd$ in case of $\eta=3,4,5$.}
\label{P1}
\end{figure}

Next, we compare our results to an equidistant sampling method, more precisely to the one of Cohen and Lindner \cite{CL:13}. Sampling at equidistant times $\Delta, 2\Delta, \dots, n\Delta$ for $\Delta > 0$, leads to an autocovariance function 
$$\gamma_{eq}(h) = \ew(X_0X_h) = \sigma^2 \int_\R f(u)f(u+h) \rmd u = \frac{\sigma^2}{2a} \e^{-ah} \, , \quad h > 0 \, ,$$
from which we conclude that $\rho_{eq}(\Delta)=\gamma_{eq}(\Delta)/\gamma(0) = \e^{-a\Delta}$ and hence
\begin{align}\label{eq5.3}
a=-\frac{\log(\rho_{eq}(\Delta))}{\Delta} \, .
\end{align}
For an estimator of $\rho_{eq}(\Delta)$, i.e. for $\rho_{eq}^*(\Delta) = \gamma_{eq;n;\Delta}^*(\Delta)/\gamma_{eq;n;\Delta}^*(0)$, where $\gamma_{eq;n;\Delta}^*(h\Delta)=\frac{1}{n} \sum_{t=1}^n X_{t\Delta}X_{(t+h)\Delta}$, $h \in \N$, by Theorem 3.5 of Cohen and Lindner \cite{CL:13} we have
\begin{align}
\sqrt{n} (\rho_{eq}^*(\Delta) - \rho_{eq}(\Delta)) \dcon N(0,V) \, , \quad n \to \infty \, ,
\end{align}
where 
\begin{align}
V &= \frac{(\eta-3)\sigma^4}{\gamma_{eq}(0)^2} \int_0^\Delta (g_{1;\Delta}(u)-\rho(\Delta)g_{0;\Delta})^2 \rmd u \nn \\
&\quad + \sum_{k=1}^\infty (\rho((k+1)\Delta) + \rho((k-1)\Delta) - 2\rho(\Delta)\rho(k\Delta))^2 \nn
\end{align}
\noindent with $g_{q;\Delta} \colon [0,\Delta] \to \R \, , \ u \mapsto \sum_{k=-\infty}^\infty f(u+k\Delta)f(u+(k+q)\Delta)$
given as in Proposition 3.1 of Cohen and Lindner \cite{CL:13}.

Knowing this, we suggest as an estimator of $a$ given in (\ref{eq5.3}) $\widehat{a}_{eq} := -\frac{\log(\rho_{eq}^*(\Delta))}{\Delta}$.
A simple calculation of $V$ for the specific kernel and the specific autocovariance function, and an application of the delta-method then leads to
$$\sqrt{n} (\widehat{a}_{eq} - a) \dcon N(0,\Delta^{-2} (\e^{2a\Delta}-1)) \, , \quad n \to \infty \, .$$

For comparing the asymptotic variances of the estimators $\widehat{a}$ and $\widehat{a}_{eq}$, we select different time scales by choosing $\Delta=\frac{1}{\lbd}$. We remind that in the renewal sampling case the expected waiting times between two sample times $T_i$ and $T_{i+1}$ is given by $\ew(W_1)=\frac{1}{\lbd}$. Then the asymptotic relative efficiency $\sigma^2_{eff}$ is given by
$$\sigma^2_{eff}=\frac{\frac{(\lbd+a)^4}{\lbd^2} \mathbf{W}_{11} - a^2}{\lbd^{2} (\e^{2a\frac{1}{\lbd}}-1)} \, ,$$
where $\mathbf{W}_{11}$ given as in (\ref{*r53}). 

We plot in Figure \ref{P1} the relative efficiency $\sigma^2_{eff}$ with respect to the mean reverting parameter $a$. The dotted line belongs to $\eta=3$, the solid line to $\eta=4$ and the dashed line to $\eta=5$.




\begin{table}[ht]
\centering
\begin{tabular}{|l|c|c|c|}
\hline
\diagbox{$\lbd$}{$\eta$} & $3$ & $4$ & $5$ \\
\hline \hline
$0.05$ & $0.1288$ & $0.1294$ & $0.1300$ \\
$0.5$ & $1.2878$ & $1.3455$ & $1.3983$ \\
$1$ & $2.5755$ & $2.7965$ & $2.9814$ \\
$2$ & $5.1509$ & $5.9627$ & $6.5465$ \\
\hline
\end{tabular}
\caption{Values of $a$ for which $\widehat{a}$ becomes more efficient depending on $\lbd$ and $\eta$.}
\label{P2}
\end{table}

The estimator $\widehat{a}$ is more efficient than $\widehat{a}_{eq}$ as $a$ tends to infinity. Table \ref{P2} shows, depending on $\lbd$ and $\eta$, the smallest value of $a$ for which $\sigma^2_{eff} \leq 1$. For values of $a$ less than $2$, the estimator based on an equidistant sampling is more efficient than $\widehat{a}$ unless the sampling frequency $\Delta$ is greater than $1$. 

We see that the non-equidistant sampling performs worse as the kurtosis of the driving L\'evy process increases. The best scenario across all time scales is observed for $\eta=3$ which corresponds to the Brownian motion case.

\section{Conclusion}
\setcounter{equation}{0}

In this paper, we studied distributional limits for the sample mean and the sample autocovariance and autocorrelation functions of a L\'evy driven continuous time moving average process. We achieved these results by assuming slightly more restrictive conditions with respect to the work of Cohen and Lindner \cite{CL:13} and gave an application of the theory in estimating the mean reverting parameter of a L\'evy driven OU process.

Of particular interest is investigating the asymptotic behavior of the sample moments when the renewal sampling is not independent of the driving L\'evy process which we hope to address in future work.

\vspace{0.2cm}
\noindent{\bf Acknowledgment}
\vspace{0.2cm}

We would like to deeply thank Alexander Lindner for many hours of discussion and his helpful comments and remarks.

\end{document}